\documentclass[leqno,11pt]{amsart}

\usepackage{amsmath}
\usepackage{amssymb}
\usepackage{enumerate,epsfig}

\newtheorem{theorem}{Theorem}[section]
\newtheorem{proposition}[theorem]{Proposition}
\newtheorem{lemma}[theorem]{Lemma}
\newtheorem{corollary}[theorem]{Corollary}

\theoremstyle{remark}
\newtheorem{remark}[theorem]{Remark}

\newcommand{\R}{ {\mathbb R} }

\newcommand{\Z}{{\mathbb Z}}
\newcommand{\T}{{\mathbb T}}
\newcommand{\cD}{{\mathcal D}}
\newcommand{\cB}{{\mathcal B}}
\makeatletter
\@addtoreset{equation}{section}
\makeatother

\newcommand{\cqfd}{{\unskip\kern 6pt\penalty 500
\raise -2pt\hbox{\vrule\vbox to 6pt{\hrule width 6pt
\vfill\hrule}\vrule}\par}}
\newcommand{\del}{\partial}
\newcommand{\dt}[1]{{\displaystyle \frac{\partial #1}{\partial t}}}
\newcommand{\ep}{{\varepsilon}}
\newcommand{\supp}{{\rm supp \,}}
\newcommand{\dist}{{\rm dist \,}}
\newcommand{\eps}{{\varepsilon}}

\newcommand{\delp}{{\delta_{p}}}
\def\charf {\mbox{{\text 1}\kern-.24em {\text l}}}
\renewcommand{\div}{\textrm{div}}

\renewcommand{\det}{\hbox{det}\,}
\newcommand{\tr}{\mathop{\mathrm{tr}}}
\numberwithin{equation}{section}

\begin{document}
\title [ kinetic homogenization]
{
Kinetic decomposition for periodic homogenization problems
}

\author[P.-E. Jabin]{Pierre-Emmanuel Jabin}
\address{Pierre-Emmanuel Jabin -- Lab. J.-A. Dieudonn\'e ,
 CNRS UMR 6621, Univ. Nice - Sophia-Antipolis, Parc Valrose, 
 06108 Nice CEDEX 02, France}
 \email{ jabin@math.unice.fr}

\author[A.E. Tzavaras]{Athanasios E. Tzavaras}
\address{Athanasios E. Tzavaras -- Department of Mathematics,
University of Maryland, College Park, Maryland 20742, USA 
}
\email{tzavaras@math.umd.edu}

\begin{abstract}
We develop an analytical tool which is adept for detecting
shapes of oscillatory functions, is useful in decomposing
homogenization problems into limit-problems for kinetic equations, and provides 
an efficient framework for the validation of multi-scale asymptotic expansions. 
We apply it first to a hyperbolic homogenization problem 
and transform it to a hyperbolic limit problem for a kinetic equation.
We establish conditions determining an effective equation
and counterexamples for the case that such conditions fail. 
Second, when the kinetic decomposition is applied to the problem of enhanced diffusion,
it leads to a diffusive limit problem for a kinetic equation that in turn
yields the effective equation of enhanced diffusion.
\end{abstract}


\maketitle
\baselineskip 16pt

%
%
%
%

\section{Introduction}
\label{intro}

Homogenization problems appear in various contexts of science and engineering and involve
the interaction of two or more oscillatory scales. In this work we focus on the simplest possible
mathematical paradigms of periodic homogenization. 
Our objective is to develop an analytical
tool that is capable of understanding the shapes of periodic oscillatory functions when the scales
of oscillations are a-priori known (or expected), and use it in order to transform the homogenization
problem into a limit problem for a kinetic equation. The calculation of an effective equation
becomes then an issue of studying a hyperbolic (or diffusive) limit for the kinetic equation. The 
procedure is well adapted in identifying the specific characteristics of the underlying homogenization
problem and provides an efficient tool for the rigorous justification of multiscale asymptotic
expansions. 

The main idea is motivated from considerations of kinetic theory. When the statistics of interacting particles
is studied it is customary to introduce an empirical measure and to study its statistical properties in the
(weak) limit when the number of particles gets large. Likewise, for an oscillating family of
functions $\{ u_{\eps} \}$ if we want to study the shape of periodic oscillations at a predetermined scale 
we may introduce an inner variable that counts the content of oscillation at such scale. 
For instance, to count oscillations at the scale $\frac{x}{\eps}$ one can introduce
\begin{equation}
\label{intro1}
f_{\eps} (x, v) = u_{\eps}(x) \delta_{p} (v - \frac{x}{\eps})
\end{equation}
where $\delta_{p}$ is the periodic delta function, and study the family $\{ f_{\eps} \}$. A-priori bounds
for $\{ u_{\eps} \}$ translate to uniform bounds for $\{ f_{\eps} \}$:  if for example
$u_{\eps}$ is uniformly bounded in $L^{2}$, $u_{\eps} \in_{b} L^{2}$, then 
$f_{\eps} \in_{b} L^{2} (M(\T^{d}))$ and, along a subsequence,
\begin{equation}
\label{intro2}
f_{\eps} \rightharpoonup f \quad \text{weak$\star$ in $L^{2} (M(\T^{d}))$} \, ,
\end{equation}
where $M(\T^{d})$ stands for the periodic measures.
In addition, the resulting $f$ is better: $f \in L^{2}(L^{2}(\T^{d}))$.

The above object should be compared to the concept of double-scale limit
introduced in the influential work of Nguetseng \cite{Nguetseng}
and applied to a variety of homogenization problems \cite{Allaire,E,HX,FP}.
In the double-scale limit one tests the family 
$\{ u_{\eps} \}$ against oscillating test functions and develops a 
representation theory for the resulting weak-limits.
It turns out, \cite{Nguetseng}, that for a uniformly bounded family 
$u_{\eps} \in_{b} L^{2}$ and test functions $\varphi$ periodic in $v$ 
\begin{equation}
\label{intro3}
\int u_{\eps} (x) \varphi(x,\frac{x}{\eps} )dx 
\to \int \int f(x,v) \varphi(x,v) dx dv
\end{equation}
where $f \in L^{2}(L^{2}(\T^{d}))$. The reader should note that 
this is precisely the content of \eqref{intro1}, \eqref{intro2},
which thus provide an alternative interpretation to the double scale limit.
However, what seems to have been missed, perhaps because
Nguetseng's analysis \cite{Nguetseng} proceeds without writing down
\eqref{intro1} but rather by establishing directly \eqref{intro3},
is that the measures $f_{\eps}$ satisfy in their own right very
interesting equations.
This is a consequence of additional properties, like
\begin{equation}
\label{intro4}
\Big ( \nabla_{x} + \frac{1}{\eps} \nabla_{v} \Big ) f_{\eps} (x, v) =
\nabla_{x} u_{\eps}(x) \delta_{p} (v - \frac{x}{\eps}) \, ,
\end{equation}
obtained by applying differential operators that annihilate 
the singular measure.
Properties like \eqref{intro4},  in turn,  suggest a procedure for embedding
homogenization problems into limit problems for kinetic equations.
In the sequel we develop this perspective,
using as paradigms the problem of
hyperbolic homogenization, and the problem of enhanced diffusion.

The double-scale limit \cite{Nguetseng} along with
the technique of multiscale asymptotic expansions  \cite{BLL} have been quite
effective in the development of homogenization theory
with considerable progress in several contexts ({\it e.g.}
\cite{MPP}, \cite{Allaire}, \cite{AM}, \cite{E},
\cite{FP}, \cite{HX}).
Other tools have also been used for the
homogenization of linear hyperbolic problems:
Among them are of course Young measures, developed by Tartar
and used for the homogenization of some particular linear transport
equations in two dimensions (see \cite{Tartar1} and
\cite{Tartar2}). Wigner measures (see \cite{GMMP}) may also be mentioned.

As our first example we consider the hyperbolic homogenization problem
\begin{equation}
\label{intro5}
\begin{aligned}
\dt {u_{\eps}} + a \big (x, \frac{x}{\eps} \big) \cdot \nabla_{x} u_{\eps} &= 0 
\\
u_{\eps} (0, x) &= U^{0}(x, \frac{x}{\eps}) \, ,
\end{aligned}
\end{equation}
with $a(x,v)$ a divergence free field periodic in $v$,  is transformed to the problem of
identifying the hyperbolic limit $\eps \to 0$ of the kinetic initial-value problem
\begin{equation}
\label{intro6}
\begin{aligned}
\dt {f_\ep} +a(x,v)\cdot\nabla_x f_\ep +\frac{1}{\ep}a(x,v)\cdot
\nabla_v f_\ep &=0,
\\
f_\ep(t=0,x,v) 
&=U^0(x,v)\;\delta_{p} (v-\frac{x}{\eps})
\end{aligned}
\end{equation}
Homogenization for \eqref{intro5} has been studied by Brenier
\cite{Brenier},
E \cite{E}, Hou and Xin \cite{HX} and, in fact,
the effective equation is sought - motivated by the double-scale limit
- in a class of kinetic equations. Eq. \eqref{intro5} is by no means
the only interesting hyperbolic problem for homogenization; we refer
to \cite{AHZ}, \cite{Hamdache},  
\cite{GouPou} (where a kinetic equation itself is homogenized), and
to
\cite{AlVa} for an example concerning a Schr\"odinger equation (the list is
of course not exhaustive).

For \eqref{intro5}, our analysis proceeds by studying
the hyperbolic limit for the kinetic equation \eqref{intro6}.
We find that if the kernel of the cell-problem
\begin{equation}
\label{introh1}
K_x=\Big\{g\in L^2(\R^d\times\T^d)\quad  \Big|\ a(x, v)\cdot
\nabla_v g=0\ in\ \cD'\Big\}
\end{equation}
is {\it independent} of $x$, then it is possible to identify 
the effective equation. 
Namely, when the vector fielfd $a = a(v)$ is independent of $x$ 
the effective equation for $f$ reads
\begin{equation}
\label{intro7}
\begin{aligned}
\dt f + (P a) \cdot \nabla_{x} f &= 0
\\
f(t=0, x,v) &= P U^{0} (x,v) \, ,
\end{aligned}
\end{equation}
where $P$ is the projection operator on the kernel $K$, and in turn
$u = \int_{\T^d} f dv$  (see Theorem \ref{hyphomthm}). 
By contrast, when $a = a(x)$ and $K_{x}$ depends on $x$, a counterexample
is constructed that shows that the effective equation can not be 
a pure transport equation (see section \ref{secce}). 
In section \ref{sectr}, this analysis is extended for homogenization
problems where a periodic fine-scale structure is transported 
by a divergence-free vector field (see equations \eqref{tos} and \eqref{kintos})
analogous results to the case of \eqref{intro5} are found. 
Such kinetic equations might turn very useful for devising computational 
algorithms for the computation of homogenization problems.

A second paradigm is the problem of enhanced diffusion
\begin{equation}
\label{intro8}
\begin{aligned}
\del_{t} u_{\eps} + \frac{1}{\eps} a(x, \frac{x}{\eps}) \cdot \nabla_{x} u_{\eps}
&= \alpha \triangle_{x} u_{\eps}
\\
u_{\eps} (0, x) = U^{0} (x, \frac{x}{\eps})
\end{aligned}
\end{equation}
with $a(x,v)$ periodic, divergence-free and with mean $\int_{\T^{d}} a = 0$. The results formally
obtained by multiscale asymptotics have been validated for this problem by 
McLaughlin, Papanicolaou and Pironneau \cite{MPP}, 
Avellaneda and Majda \cite{AM},
and Fannjiang and Papanicolaou \cite{FP}. 
We revisit this problem from the perspective of the
kinetic decomposition and transform it to the problem of 
identifying the $\eps \to 0$ limit 
\begin{equation}
\label{intro9}
\begin{aligned}
\dt {f_\ep} + \frac{1}{\eps} a(x, v)\cdot\nabla_x f_\ep 
&+ \frac{1}{\ep^{2}} 
\big ( a(x,v) \cdot \nabla_{v} f_{\eps} - \alpha \triangle_{v} f_{\eps} \big )
\\
 &= \alpha \triangle_{x} f_{\eps} 
 + \frac{2\alpha}{\eps} \nabla_{x} \cdot \nabla_{v} f_{\eps}
\, , 
\\
f_\ep(t=0,x,v) 
&=U^0(x,v)\;\delta_{p} (v-\frac{x}{\eps}).
\end{aligned}
\end{equation}
The latter is a limit for the transport-diffusion equation \eqref{intro9}
in the so-called diffusive scale, and its analysis provides the 
effective equation \eqref{basic}-\eqref{cellp}
of enhanced diffusion (see Theorem \ref{enhdiffthm}). This example indicates the efficiency of 
this approach in the rigorous validation of multi-scale asymptotic expansions.

Finally, we note that the scales of the drift and of the diffusion 
in \eqref{intro8} may be chosen differently from $1/\eps$ and $\alpha$, 
yielding other interesting homogenization problems, see for instance
Capdeboscq \cite{Capdeboscq1,Capdeboscq2}. 

The article is organized as follows. Analytical considerations like the proper definition of \eqref{intro1},
the characterization of the weak limit points of $f_{\eps}$ under various uniform bounds, 
the differential relations such as \eqref{intro4}, and the identification of asymptotics for $f_{\eps}$ 
are developed in section \ref{msde} and in appendix I. In section \ref{homohyp}, we study the
hyperbolic homogenization problem \eqref{intro5},  derive the effective equation, and
produce the counterexample mentioned before. 
Some material from ergodic theory needed in the derivation is outlined in the appendix II.
In section \ref{sectr}, we study the transport via a divergence-free field depending
on an oscillating fine-scale, we derive the associated kinetic equation, and discuss
the connection of the two formulations via characteristics and the derivation of an
effective equation. Finally, in section \ref{endiff} we study the parabolic homogenization
problem \eqref{intro8} and derive the enhanced diffusion equation via the kinetic decomposition.

%

\section{Multi-scale decomposition}
\label{msde}

Let $\{ u_{\eps}(x) \}$ be a family of functions defined on a open set 
$\Omega \subset \R^{d}$
that contains periodic oscillations and suppose that the scales of oscillations 
are either a-priori known (or anticipated).
Our goal is to introduce an analytical object that 
will prompt the anticipated scale(s) of oscillations and quantify the structure 
of oscillations in the family at the preselected scale(s).

Suppose that  periodic oscillations of length $\eps$ are anticipated in 
the family $\{ u_{\eps} \}$. 
To focus on them we consider a periodic grid with sides of length $\eps$
in each coordinate direction. The grid splits the Euclidean space into distinct cubic cells 
of volume $\eps^{d}$, and it is arranged so that the centers of the cells 
occupy the lattice $\eps \Z^{d}$. Let $\Omega$ be placed on that grid,
and define a function $\chi_{\eps} : \Omega \to \eps \Z^{d}$ that maps the generic 
$x \in \Omega$ to the center $\chi_{\eps}(x)$ of the cell containing $x$. 
To each point $x \in \Omega$ there is associated a decomposition $(\chi_{\eps}(x), v)$ 
where $\chi_{\eps}(x) \in \eps \Z^{d}$ stands for the center of the cell that $x$ occupies,
and $v \in \T^{d}$ is the vector difference $x-\chi_{\eps}(x) $  as measured in units 
of distance $\eps$,  that is  $x = \chi_{\eps}(x) + \eps v$. 
We introduce the quantity
\begin{equation}
\label{decompose}
f_{\eps}(x, v) = u_{\eps} (x) \delta_{p} \big (v - \frac{x - \chi_{\eps}(x)}{\eps} \big )
\, , \quad x \in \R^{d} \, , \; v \in \T^{d} \, ,
\end{equation}
where $\delta_{p}$ stands for a periodization of the usual delta function
with period $1$ in each coordinate direction, and $\T^{d}$ stands for the d-dimensional torus,
the quotient of $\R^{d}$ by the subgroup $\Z^{d}$.

We note that the map $x \mapsto (\chi_{\eps}(x), v) $ is single valued for points that
fall into a single cell, but multi-valued for points that fall onto the boundaries between
adjacent cells. For the latter points there would be two different decompositions 
$(\chi_{\eps}, v)$ and $(\chi'_{\eps}, v')$ associated to the same point $x \in \Omega$.
Nevertheless, in that case $x = \chi_{\eps} + \eps v = \chi'_{\eps} + \eps v'$ and,
due to the use of a periodic delta function,
$$
\delta_{p} \big (v - \frac{x - \chi_{\eps}}{\eps} \big ) = 
\delta_{p} \big (v - \frac{x }{\eps} \big ) =
\delta_{p} \big (v - \frac{x - \chi'_{\eps}}{\eps} \big )
$$
Hence, both decompositions provide 
the same outcome in \eqref{decompose} with $f_{\eps}$ defined
for $x \in \Omega$ and $v \in \T^{d}$.

The operator $\nabla_{x} + \frac{1}{\eps} \nabla_{v}$ annihilates
the form $v - \frac{x }{\eps}$ and that - at least formally - yields the formula
\begin{equation}
\label{formder}
\big ( \nabla_{x} + \frac{1}{\eps} \nabla_{v} \big ) f_{\eps} =
(\nabla_{x} u_{\eps}) \delta_{p} \big (v - \frac{x }{\eps} \big )
\end{equation}
In the sequel, we provide formal definitions
for the decomposition \eqref{decompose} and extensions as well as
differentiation properties like \eqref{formder} 
that are helpful in later sections for validating multiscale expansions.

\subsection{Definitions}
We make extensive use of distributions defined on the torus $\T^{d}$. 
Such distributions are in one-to-one correspondence  with periodic 
distributions 
$T$ on $\R^{d}$ of period 1 in each coordinate direction,
that is distributions satisfying for $i = (i_{1}, ... , i_{d}) \in
\Z^{d}$ the 
property
$\tau_{i}T = T$ where $\tau_{i}$ is the shift operator, 
see \cite[p. 229]{Schwartz}. 
The same notation is used for both interpretations of periodic distributions.
Let  $\delta_{p}$ be the periodic delta function of period 1, defined
by its action $< \delta_{p} , \psi > = \psi(0)$
on continuous periodic test functions $\psi \in C(\T^{d})$. 

We use the notation $C_{p} = C (\T^{d})$ for the continuous periodic functions, 
$C^{\infty}_{p}= C^{\infty}(\T^{d})$ for periodic test functions
and $M_{p} = M^{1} (\T^{d})$ for the periodic measures, 
with period $1$ in each coordinate direction.
Recall that $C_{p}$ is separable and that bounded sets in
$M_{p} =  \big (C_{p} \big )^{*}$ are sequentially precompact in the weak-$\star$
topology of $M_{p}$.

\subsubsection{The double-scale kinetic decomposition}
\label{ss11}
Our objective is to define the product
\begin{equation}
\label{decomp}
f_{\eps}(x, v) = u_{\eps} (x) \delta_{p} \big (v - \frac{x }{\eps} \big )
\end{equation}
which, in view of the periodicity of $\delta_{p}$ and  
$\frac{\chi_{\eps}}{\eps} \in \Z^{d}$,
coincides with \eqref{decompose}. Products of distributions are not in general well defined.
However, due to the tensor-product-like structure the product in \eqref{decomp} is well
defined by using the Schwartz kernel theorem \cite[Ch V]{Hormander}. We defer the details
for the example in section \ref{ss12}, and note that 
the definition of \eqref{decomp} 
is effected for $u_{\eps} \in L^{1}_{loc}(\Omega)$ by setting
\begin{equation}
<f_{\eps}, \theta> = \int_{\R^{d}} u_{\eps}(x) 
     \theta \big (x , \frac{x }{\eps} \big ) \, dx
\end{equation}
for $\theta (x,v) \in  C^{\infty}_{c} (\Omega ; C^{\infty} { (\T^{d})})$.
Moreover, we have the differentiation formula:

\begin{lemma}
\label{lem1}
For  $u_{\eps} \in W^{1,1}_{loc}(\Omega)$,
$$
\big ( \nabla_{x} + \frac{1}{\eps} \nabla_{v} \big ) 
\left (u_{\eps}  \delta_{p} \big (v - \frac{x }{\eps} \big )  \right )  =
(\nabla_{x} u_{\eps}) \delta_{p} \big (v - \frac{x }{\eps} \big )
$$
\end{lemma}

\begin{proof}
For $\theta (x,v) \in  C^{\infty}_{c} (\Omega ; C^{\infty}_{p})$, we 
have
$$
\begin{aligned}
< \big ( \nabla_{x} + \frac{1}{\eps} \nabla_{v} \big )  f_{\eps} , \theta>
&=
- < f_{\eps} , \big ( \nabla_{x} + \frac{1}{\eps} \nabla_{v} \big ) \theta>
\\
&= - \int_{\R^{d}} u_{\eps}(x) 
    \Big (\nabla_{x} \theta  + \frac{1}{\eps} \nabla_{v} \theta \Big)
     \big (x , \frac{x }{\eps} \big ) \, dx
\\
&= - \int_{\R^{d}} u_{\eps}(x)  \nabla_{x} 
\big (  \theta \big (x , \frac{x }{\eps} \big ) \big ) \, dx
\\
&= \int_{\R^{d}} (\nabla_{x} u_{\eps} )(x)  
 \theta \big (x , \frac{x }{\eps} \big )  \, dx
\\
&=  < (\nabla_{x} u_{\eps} ) \delta_{p} \big (v - \frac{x }{\eps} \big ) , \theta >
\end{aligned}
$$
\end{proof}

\subsubsection{A generalization}
\label{ss12}
Let  $\alpha (x) $ be a smooth vector field and
$u_{\eps} \in L^{1}_{loc}(\Omega)$. We proceed to define the product 
\begin{equation}
\label{vardecomp}
f_{\eps} =  u_{\eps} (x) \delta_{p} \big (v - \frac{\alpha (x) }{\eps} \big ) \, .
\end{equation}
Naturally it should act on tensor products $\varphi \otimes \psi$ of test functions via
the formula
\begin{equation}
\label{deftensor}
< f_{\eps}, \varphi \otimes \psi > =
\int_{\R^{d}} u_{\eps}(x) \varphi (x) \psi \big( \frac{\alpha (x) }{\eps} \big )
\, dx \, .
\end{equation}
To define \eqref{vardecomp}, we employ the Schwartz kernel theorem 
\cite[Thm 5.2.1]{Hormander}. 
Consider the linear map 
$$
\mathcal{K} : C^{\infty}(\T^{d}) \to \mathcal{D}'(\Omega) 
\quad \text{defined by } \quad \mathcal{K} \psi = u_{\eps}(x) 
\psi \big ( \frac{\alpha(x)}{\eps} \big )
$$
If $\psi_{n} \to 0$ in $C^{\infty}(\T^{d})$ then $\mathcal{K} \psi_{n} \to 0$ 
in $\mathcal{D}'(\Omega)$. The kernel theorem implies that there exists a {\it unique} 
distribution $K$ such that $< K , \varphi \otimes \psi > = ( \mathcal{K} \psi , \varphi)$,
that is $K$ acts on tensor products via \eqref{deftensor} and is the desired product. 
It satisfies, for $\theta \in C^{\infty}_{c}(\Omega ; C^{\infty}(\T^{d}))$,
\begin{equation}
\label{defn2}
< u_{\eps} \delta_{p} \big (v - \frac{\alpha (x) }{\eps} \big ), \theta > =
\int_{\R^{d}} u_{\eps}(x) \theta \big(x,  \frac{\alpha (x) }{\eps} \big ) \, dx
\, ,
\end{equation}
which can also serve as a direct definition of $f_{\eps}$.
Of course smoothness of $a(x)$ is required for the above definition: at
least $a \in C(\Omega ; \R^d)$ if $f_\eps$ is interpreted as a measure, and
more if $f_\eps$ is interpreted as a distribution and we need to take derivatives.

We now prove.

\begin{lemma}
\label{lem2}
Let  $u_{\eps} \in W^{1,1}_{loc}(\Omega)$ and
 $\alpha \in C^1 (\Omega ; \T^{d})$. Then
$$
\big ( \nabla_{x} + \frac{1}{\eps} (\nabla \alpha)^{T} \nabla_{v} \big ) 
\left (u_{\eps}  \delta_{p} \big (v - \frac{\alpha(x)}{\eps} \big )  \right )  =
(\nabla_{x} u_{\eps}) \delta_{p} \big (v - \frac{\alpha(x) }{\eps} \big )
$$
\end{lemma}

\begin{proof}
For the $k$-th coordinate, we have
$$
\begin{aligned}
< 
\Big ( \del_{x_{k}} 
&+ \frac{1}{\eps} \sum_{j} \frac{\del \alpha_{j}}{\del x_{k}} \del_{v_{j}}
\Big ) u_{\eps} \delta_{p} \big (v -  \frac{\alpha(x)}{\eps} \big ) , \theta 
>
\\
&=
- < u_{\eps} \delta_{p} \big (v - \frac{\alpha(x)}{\eps} \big ) , 
\del_{x_{k}} \theta  + \frac{1}{\eps} \sum_{j} \del_{v_{j}} 
\Big ( \frac{\del \alpha_{j}}{\del x_{k}} \theta \Big )
>
\\
&= - \int_{\R^{d}} u_{\eps} (x)
\del_{x_{k}} \Big ( \theta(x, \frac{\alpha(x)}{\eps} ) \Big ) \, dx
\\
&= \int_{\R^{d}} (\del_{x_{k}}u_{\eps}) (x)
\ \theta(x, \frac{\alpha(x)}{\eps} ) \, dx
\\
&= < (\del_{x_{k}}u_{\eps}) \delta_{p} \big (v - \frac{\alpha(x)}{\eps} \big )
, \theta>
\end{aligned}
$$
\end{proof}

\subsubsection{A multiscale kinetic decomposition}
\label{ss13}
We pursue next the construction of decompositions in cases when more than 
two scales are involved. Suppose that for an oscillating family $\{ u_{\eps} \}$ we wish to 
focus on oscillations at the scales $1$, $\frac{x}{\eps}$ and $\frac{x}{\eps^{2}}$. 
We define
\begin{equation}
\label{defn3sc}
f_{\eps} (x,v,w) = 
u_{\eps} (x) 
\delta_{p} \big( v - \frac{x }{\eps} \big)
\delta_{p} \big( w - \frac{v }{\eps} \big) \, , 
\qquad x \in \Omega, \, v\in \T^{d}, w\in\T^{d} \, ,
\end{equation}
or, in terms of the action on test functions,
$\theta (x,v,w) \in
C^{\infty}_{c}(\Omega ; C^{\infty} (\T^{d}\times \T^{d}))$ via the formula
\begin{equation}
\label{defnpro3sc}
<f_{\eps} , \theta> = \int_{\R^{d}} 
\theta \big(x, \frac{x}{\eps} , \frac{x}{\eps^{2}}\big) \, dx
\, .
\end{equation}
In a straightforward generalization of Lemma \ref{lem1}, $f_{\eps}$ satisfies, for
$u_{\eps} \in W^{1,1}_{loc}(\Omega)$, the
differentiation formula
\begin{equation}
\begin{aligned}
\big (  \nabla_{x} + \frac{1}{\eps} \nabla_{v} 
+ \frac{1}{\eps^{2}} \nabla_{w} \big)
f_{\eps} 
=
(\nabla_{x} u_{\eps })(x)
\delta_{p} \big( v - \frac{x}{\eps} \big)
\delta_{p} \big( w - \frac{v}{\eps} \big)
\, .
\end{aligned}
\end{equation}

To motivate the definition \eqref{defn3sc} consider for simplicity the
case that $1/\eps$ is an integer. Fix a first grid of size $\eps$ and introduce 
the quantities $\chi_{\eps}(x)$ and  
$v = \frac{x - \chi_{\eps}}{\eps} \in \T^{d}$ as before.
To focus on the scale $\frac{x}{\eps^{2}}$ we consider a second grid of length $\eps^{2}$
embedded in the first grid. When $1/ \eps$ is an integer, 
the grids fit perfectly onto one-another. 
Define the function $\psi_{\eps} : \T^{d} \to \eps \Z^{d}$ that takes the generic point $v$ 
to the center of the inner cell containing $v$, and introduce a second inner variable 
$w = \frac{v - \psi_{\eps}(v)}{\eps} \in \T^{d}$ describing the vector distance between
$v$ and the center of the inner cell containing $v$ in units of length $\eps$.  The process defines
a decomposition of the physical space $x \mapsto (\chi_{\eps}(x), v, \psi_{\eps}(v), w)$, 
and allows to define a kinetic function representing three scales by
$$
\begin{aligned}
f_{\eps} (x,v,w) &= u_{\eps} (x) 
\delta_{p} \big( v - \frac{x - \chi_{\eps}}{\eps} \big)
\delta_{p} \big( w - \frac{v - \psi_{\eps}}{\eps} \big)
\\
&=
u_{\eps} (x) 
\delta_{p} \big( v - \frac{x }{\eps} \big)
\delta_{p} \big( w - \frac{v }{\eps} \big) \, , 
\qquad x \in \Omega, \, v\in \T^{d}, w\in\T^{d} \, ,
\end{aligned}
$$
This definition is also good when $1/ \eps$ is not an integer as can be seen
by the formula \eqref{techn2} in the appendix.

\subsection{Multiscale analysis of uniformly bounded families of functions}
Nguetseng \cite{Nguetseng} introduced the notion of double scale limit,
which has been a very effective technical tool in the development of periodic homogenization
theory. 
His approach does not use the kinetic decomposition \eqref{decompose}, but
the double-scale limit is precisely the weak limit of the measures introduced 
in \eqref{decompose}.
For this reason, we will review the results of Nguetseng \cite{Nguetseng}
from the perspective of the theory presented here, and produce some further asymptotic analysis
of kinetic decompositions for uniformly bounded families of functions.
In the sequel, the notation $u_{\eps} \in_{b} X$ means that the family $\{u_{\eps}\}$
belongs in a bounded set of the Banach space $X$.

\subsubsection{Uniform $L^{2}$-bounds}

Suppose first that $\{u_{\eps}\}$ satisfies $u_{\eps} \in_{b} L^{2}(\Omega)$. 
We define $f_{\eps}$ by \eqref{decomp} and note that 
\begin{equation}
\label{l2bound}
f_{\eps} \in_{b} L^{2} (\Omega  ; M_{p}) \, .
\end{equation}
The Riesz representation theory asserts that there is an isometric isomorhism between the
dual of $C_{p}=C(\T^{d})$ and the Banach space of periodic Radon measures 
$M_{p} = M^{1} (\T^{d})$ on the torus. Since $C_{p}$ is separable, bounded sets in $M_{p}$
are sequentially precompact in the weak-$\star$ topology of $M_{p}$. Also, since $C_{p}$
is separable, so is $L^{2} (\Omega ; C_{p} )$ and thus bounded sets in $L^{2} (\Omega ; M_{p})$
are sequentially precompact in the weak-$\star$ topology of $L^{2} (\Omega ; M_{p})$. 

As a consequence \eqref{l2bound} implies that, along a subsequence,
\begin{equation}
\label{weakst}
f_{\eps} \rightharpoonup f  \quad \text{ weak-$\star$ in $L^{2} (\Omega ; M_{p})$}
\end{equation}
with $f \in L^{2} (\Omega ; M_{p})$, that is 
\begin{equation}
\label{weakstdet}
\begin{aligned}
< f_{\eps} , \theta > &= \int u_{\eps}(x) \theta (x, \frac{x}{\eps}) \, dx
\\
&\to <f, \theta> = \iint f(x,v) \theta(x,v) \, dx dv  
\qquad \text{ for $\theta \in L^{2} (\Omega ; C_{p})$}
\end{aligned}
\end{equation}

\bigskip
\noindent
{\sc Examples.}
A few examples will illustrate the properties of this convergence.

1. 
Note first that
\begin{equation}
\label{ex1}
\delta_{p} \big( v - \frac{x}{\eps} \big) \rightharpoonup 1 
\qquad \text{weak-$\star$ in $L^{\infty}(\Omega ; M_{p})$} \, ,
\end{equation}
that is 
\begin{equation}
\label{ex1det}
< \delta_{p} \big( v - \frac{x}{\eps} \big) , \theta > =
\int_{\Omega} \theta (x, \frac{x}{\eps}) dx \to 
\int_{\Omega} \int_{\T^{d}} \theta (x,v) dx dv \, ,
\end{equation}
for $\theta \in L^{1}(\Omega ; C_{p})$. This is a classical result, see \cite{BLL},
and a proof is provided for completeness in the appendix.

2. If $u_{\eps} \to u$ strongly in $L^{2}(\Omega)$, then
\begin{equation}
\label{ex3}
u_{\eps} (x) \delta_{p} \big( v - \frac{x}{\eps} \big) \rightharpoonup u(x)
\qquad \text{weak-$\star$ in $L^{2} (\Omega ; M_{p})$}.
\end{equation}
Indeed, since $u_{\eps} \to u$ in $L^{2}(\Omega)$ and
$\psi \big( \frac{x}{\eps} \big)  \rightharpoonup \int_{T^{d}} \psi (v) dv$
weakly in $L^{2} (\Omega)$, for $\theta = \varphi \otimes \psi$ a tensor
product
$$
\begin{aligned}
< u_{\eps} \delta_{p} \big( v - \frac{x}{\eps}\big ) , \varphi \otimes \psi >
&=
\int_{\R^{d}} u_{\eps} (x)  \psi \big( \frac{x}{\eps} \big) \varphi(x) dx
\\
&\to \int_{\R^{d}} \int_{\T^{d}} u(x) \varphi(x) \psi(v) dx dv
\end{aligned}
$$
Finite sums of tensor products $\sum_{j} \varphi_{j} \otimes \psi_{j}$ are dense 
in $L^{2}(\Omega  ; C_{p} )$ and 
\eqref{ex3} follows by a density argument.

3. For  $u_{\eps} = a(\frac{x}{\eps})$. where $a(v)$ is a periodic function, we obtain
\begin{equation}
\label{weakex4}
a(\frac{x}{\eps}) \rightharpoonup \int_{\T^{d}} a(v) dv
\end{equation}
and, by \eqref{techn1},
\begin{equation}
\label{ex4}
a(\frac{x}{\eps}) \delta_{p} \big( v - \frac{x}{\eps}\big ) \rightharpoonup
a(v)
\end{equation}
Observe that the weak limit \eqref{weakex4} retains only the information of the average of $a$ 
while the double scale kinetic limit \eqref{ex4} also retains the information of the shape of $a$. 
Equation \eqref{techn2} in the appendix
indicates that if we were to prompt $u_{\eps} = a(\frac{x}{\eps})$ with test functions oscillating 
on a scale different than $\eps$ then the information on the shape of $a$ is lost and only the 
average is perceived. 
Therefore, the double scale decomposition works well when the scale of oscillations
are {\it a-priori} known, so that the right oscillatory scale is prompted.

4. Finally, let $u_{\eps} = a(\frac{x}{\eps}) b(\frac{x}{\eps^{2}})$ where 
$a$ and $b$ are periodic functions. Then \eqref{techn1} and \eqref{techn2} give
\begin{align*}
a(\frac{x}{\eps}) b(\frac{x}{\eps^{2}}) &\rightharpoonup \int_{\T^{d}} a(y)dy \;
\int_{\T^{d}} b(z)dz
\\
a(\frac{x}{\eps}) b(\frac{x}{\eps^{2}})  \delta_{p} \big( v - \frac{x}{\eps}\big )
&\rightharpoonup  a(v) \int_{\T^{d}} b(z)dz
\\
a(\frac{x}{\eps}) b(\frac{x}{\eps^{2}})  \delta_{p} \big( v - \frac{x}{\eps}\big )
\delta_{p} \big( w - \frac{v}{\eps}\big )
&\rightharpoonup  a(v) \, b(w)
\end{align*}

\bigskip
In the following proposition, we give a simplified proof of \cite[Thm 1]{Nguetseng} 
concerning the double-scale limit for sequences that are uniformly bounded in $L^{2}(\Omega)$.

\begin{proposition}
\label{prop1}
Let $u_{\eps} \in_{b} L^{2}(\Omega)$. Then, along a subsequence,
$$ 
f_{\eps} \rightharpoonup f   \qquad \text{ weak-$\star$ in $L^{2} (\Omega \, ; M_{p})$}
$$ 
with $f \in L^{2} (\Omega  \times \T^{d})$.
\end{proposition}

\begin{proof}
Let $\theta \in C^{\infty}_{c} (\Omega ; C^{\infty}_{p} )$ be a test function. Then
$$
<f_{\eps} , \theta> = \int_{\Omega} u_{\eps} (x) \theta(x, \frac{x}{\eps}) dx
$$
and 
$$
\begin{aligned}
|<f_{\eps} , \theta>|  &\le \| u_{\eps} \|_{L^{2}(\Omega)}
        \left ( \int_{\Omega} |\theta (x, \frac{x}{\eps})|^{2} dx \right )^{1/2}
\\
&\le C \left ( \int_{\Omega} |\theta (x, \frac{x}{\eps})|^{2} dx \right )^{1/2}
\\
&\stackrel{\eqref{ex3}}{\to}  C \left ( \int_{\Omega} \int_{\T^{d}} 
     |\theta (x, v)|^{2} dx dv \right )^{1/2}
\end{aligned}
$$
Hence, $f_{\eps} \in_{b} \left (  L^{2} (\Omega ; C_{p}) \right )^{*}$,
$f_{\eps} \rightharpoonup f$ weak-$\star$ in $L^{2} (\Omega ; M_{p})$
and $f \in L^{2} (\Omega ; M_{p})$. Moreover,
$$
\frac{  |<f , \theta>|  }{ \| \theta \|_{  L^{2} (\Omega \times \T^{d})  } }
=
\lim \frac{  |<f_{\eps} , \theta>|  }{ \| \theta \|_{  L^{2} (\Omega \times \T^{d})  } }
\le C
$$
and $f \in L^{2} (\Omega \times \T^{d})$.
\end{proof}

\subsubsection{ Uniform $H^{1}$-bounds}

Next consider the case of families $\{u_{\eps}\}$ that are uniformly bounded 
in $H^{1} (\Omega)$. The first proposition is essentially 
a rephrasing of \cite[Thm 3]{Nguetseng}.

\begin{proposition}
\label{prop2}
Let $u_{\eps} \in_{b} H^{1} (\Omega)$. Then, there exist $u \in H^{1} (\Omega)$,
$\pi \in L^{2} (\Omega ; H^{1}(\T^{d}))$ such that, along a subsequence,
$$
\begin{aligned}
f_{\eps} = u_{\eps} \delta_{p} \big( v - \frac{x}{\eps}  \big)
&\rightharpoonup u(x) 
\qquad  \text{  weak-$\star$ in $L^{2} (\Omega ; M_{p} )$ }
\\
(\nabla_{x} u_{\eps} - \nabla_{x} u ) \delta_{p} \big( v - \frac{x}{\eps}  \big)
&\rightharpoonup \nabla_{v} \pi (x,v) 
\qquad  \text{  weak-$\star$ in $L^{2} (\Omega ; M_{p} )$ }
\end{aligned}
$$
\end{proposition}

\begin{proof}
Along subsequences (whenever necessary) 
$u_{\eps} \rightharpoonup u$ weakly in $H^{1}(\Omega)$,
$u_{\eps} \to u$ in $L^{2}(\Omega)$ and 
$u_{\eps} \delta_{p} \big( v - \frac{x}{\eps}  \big) \rightharpoonup u$
weak-$\star$ in $L^{2}(\Omega ; M_{p} )$. 
Moreover, Proposition \ref{prop1} implies
$$
g_{i, \eps} : = 
\big(  \frac{\del u_{\eps}}{\del x_{i}} - \frac{\del u}{\del x_{i}}  \big)  
\delta_{p} \big( v - \frac{x}{\eps}  \big)
\rightharpoonup g_{i} (x,v)  
$$
weak-$\star$ in $L^{2}(\Omega ; M_{p})$ with $g_{i} \in L^{2} (\Omega \times \T^{d}) $.

For $\varphi \in C_{c}^{\infty}(\Omega)$, $\psi_{i} \in C^{\infty}_{p}$, we have
$$
\begin{aligned}
- \int_{\Omega} (u_{\eps} - u)  
  &\left [  
  \varphi (x) \frac{1}{\eps} \frac{\del \psi_{i}}{\del v_{i}} \big( \frac{x}{\eps}  \big) 
  + \frac{\del \varphi}{\del x_{i}}(x) \psi_{i} \big( \frac{x}{\eps} \big)
   \right ] \, dx dv
\\
&= \int_{\Omega} \frac{\del(u_{\eps} - u) }{\del x_{i}}
  \varphi (x) \psi_{i} \big( \frac{x}{\eps}  \big)  \; dx dv
 \\
&\to
\int_{\R^{d}} \int_{\T^{d}} g_{i}(x,v) \varphi (x) \psi_{i} (v) dx dv
\end{aligned}
$$
We apply the above formula to a test function $\Psi = (\psi_{1}, ... , \psi_{d})$ 
that satisfies $\div_{v} \Psi = 0$. 
Then
$$
- \int_{\R^{d}} (u_{\eps} - u)  
  \left [  
  \varphi (x) \frac{1}{\eps} \frac{\del \psi_{i}}{\del v_{i}} \big( \frac{x}{\eps}  \big) 
  + \frac{\del \varphi}{\del x_{i}}(x) \psi_{i} \big( \frac{x}{\eps} \big)
   \right ] \, dx dv \to 0
$$
and we conclude that for a.e. $x \in \Omega$
\begin{equation}
\int_{\T^{d}} \sum_{i} g_{i}(x,v)  \psi_{i}(v) \, dv = 0
\quad \text{for any $\Psi$ with $\div_{v} \Psi = 0$}.
\end{equation}
A lemma from \cite[Lemma 4]{Nguetseng} then implies 
there exists $\pi \in L^{2}(\Omega ; H^{1}(\T^{d}))$ such that 
$G = (g_{1}, ... , g_{d}) = \nabla_{v} \pi $.
\end{proof}

The next proposition is novel and establishes the asymptotics of $f_{\eps}$ 
when the family $\{ u_{\eps}\}$ is uniformly bounded in $H^{1}(\Omega)$.

\begin{proposition}
\label{prop3}
Let $u_{\eps} \in_{b} H^{1} (\Omega)$. Then
\begin{equation}
\label{propasym}
\begin{aligned}
\frac{1}{\eps}(\delta_{p}(v - \frac{x}{\eps}) -  1)
&\to  0 \quad \text{in $\cD'$}
\\
g_{\eps} := u_{\eps}(x) \frac{1}{\eps}(\delta_{p}(v - \frac{x}{\eps}) -  1)
 &\in_{b}  H^{-1} (\Omega ; M_{p} )
 \\
 g_{\eps} &\rightharpoonup g \quad 
 \text{ weak-$\star$ in $H^{-1} (\Omega ; M_{p})$}
 \\
 g &\in L^{2} (\Omega ; H^{1} (\T^{d})) \, ,
\end{aligned}
\end{equation}
and $f_{\eps}$ enjoys the asymptotic expansion 
$$
f_{\eps} = u_{\eps} + \eps g + o(\eps) \quad \text{in $\cD'$} \, .
$$
\end{proposition}

\begin{proof}
It is instructive to first give a quick proof for the case of one space dimension. Consider
the function $H(v) = v$, $v \in [0,1]$ and let $H_{p}$ denote its periodic extension
of period $1$.  $H_{p} \in L^{\infty}(\R)$ satisfies
 $\del_{v} H_{p} (v) = 1 - \delta_{p} (v)$ and
\begin{equation}
\label{diden}
 \frac{1}{\eps} \big ( \delta_{p} (v - \frac{x}{\eps}) - 1   \big )
=
\del_{x} H_{p} ( v - \frac{x}{\eps} )
\end{equation}
Using standard properties of weak
convergence we obtain for $\theta \in C^{1}_{c}(\Omega , C_{p})$
$$
< H_{p} (v - \frac{x}{\eps}) , \theta > 
\to
\int_{\Omega}\int_{\T^{d}} \int_{\T^{d}} H_{p} (v-w) \theta(x,v) dw dv dx
$$
and
$$
\begin{aligned}
< \frac{1}{\eps} \big ( \delta_{p} (v - \frac{x}{\eps}) - 1   \big ) , \theta >
&=
< \del_{x} H_{p} ( v - \frac{x}{\eps} ) , \theta >
\\
&=
- < H_{p} ( v - \frac{x}{\eps} ) , \del_x \theta >
\\
&\rightharpoonup
- \int_{\Omega} \int_{\T^{d}} \int_{\T^{d}} H_{p} (v-w) \del_{x} \theta(x, v) dw dv dx  
\\
&= 0 \, .
\end{aligned}
$$
Using lemma \ref{lem1} and  \eqref{diden}
we have, for $u_{\eps} \in H^{1}(\Omega)$, the identities
$$
\big ( \del_{x} + \frac{1}{\eps} \del_{v} \big ) 
\big ( u_{\eps} H_{p} (v - \frac{x}{\eps})  \big )
= (\del_{x} u_{\eps}) H_{p} (v - \frac{x}{\eps})  
$$
and
$$
u_{\eps} \frac{1}{\eps} \big ( \delta_{p} (v - \frac{x}{\eps}) - 1   \big ) 
= 
\del_{x} \Big ( u_{\eps} H_{p} (v - \frac{x}{\eps}) \Big )
- (\del_{x} u_{\eps}) H_{p} (v - \frac{x}{\eps}) \, .
$$
The remaining three properties are direct consequences of the last formula.

Consider next the multi-dimensional case.
For $\theta \in C^{1}_{c}(\Omega ; C_{p} )$, 
let $\supp \theta$ denote the support (in $x$) of  $\theta$, and fix
$\eps < \frac{1}{\sqrt{d}} \, \dist ( \supp \theta , \partial \Omega)$. 
We cover $\supp \theta$ by cubes $C_{k}$ centered at points $\chi_{k} \in \eps \Z^{d}$ 
of latteral size $\eps$.
The number of the cubes covering $\supp \theta$ is of the order 
$\eps^{-d} O(|\supp \theta|)$, and the covering is arranged so that
$\supp \theta \subset \cup_{k=1}^{N} C_{k} \subset \Omega$. We observe that
by construction  $\frac{1}{\eps} \chi_{k} \in \Z^{d}$ and compute
$$
\begin{aligned}
<g_{\eps} & , \theta >
= \int_{\Omega} u_{\eps} (x) \frac{1}{\eps}
\Big (  
\theta(x, \frac{x}{\eps}) - \int_{\T^{d}} \theta(x,v) dv
\Big )
dx
\\
&= \sum_{k \in \Z^{d}} \int_{C_{k}} u_{\eps}(x) \frac{1}{\eps}
\Big (  
\theta(x, \frac{x}{\eps}) - \int_{\T^{d}} \theta(x,v) dv
\Big )
dx
\\
&= \sum_{k \in \Z^{d}} \eps^{d} \int_{\T^{d}}
u_{\eps} (\chi_{k} + \eps w)
\frac{1}{\eps} 
\Big (  
\theta(\chi_{k} + \eps w , w) - \int_{\T^{d}} \theta(\chi_{k} + \eps w , v) dv
\Big ) dw
\\
&= \sum_{k \in \Z^{d}} \eps^{d} \int_{\T^{d}}
 \frac{1}{\eps} 
\Big (  
(u_{\eps} \theta)(\chi_{k} + \eps w , w) 
- \int_{\T^{d}} (u_{\eps} \theta)(\chi_{k} + \eps \rho , w) d\rho
\Big )
dw
\end{aligned}
$$
We employ the Poincar\'e inequality
$$
\left | v(z) - \int_{\T^{d}} v(z') dz' \right | \le
\int_{\T^{d}} | \nabla v| (z') dz'
$$
for $v(z) = (u_{\eps} \theta)(\chi_{k} + \eps z , v)$ to obtain
$$
\begin{aligned}
\frac{1}{\eps}  \Big | 
(u_{\eps} \theta)(\chi_{k} + \eps z , v) 
&- \int_{\T^{d}} (u_{\eps} \theta)(\chi_{k} + \eps \rho , v) d\rho
\Big | 
\\
&\le
\int_{\T^{d}} |\nabla_{x} (u_{\eps} \theta)| (\chi_{k} +\eps \rho , v) d\rho
\end{aligned}
$$
and
$$
\begin{aligned}
| <g_{\eps} , \theta > |
&\le 
\sum_{k \in \Z^{d}} \eps^{d} \int_{\T^{d}} \int_{\T^{d}}
| \nabla_{x} (u_{\eps} \theta) | (\chi_{k} + \eps \rho , w) \, d\rho dw
\\
&= \int_{\Omega} \int_{\T^{d}} |\nabla_{x} (u_{\eps} \theta)|(x,v) \, dxdv
\\
&\le \| u_{\eps} \|_{H^{1}(\Omega)} 
\left ( \|\theta \|_{L^{2}(\Omega \times \T^{d})} 
+ \|\nabla_{x} \theta \|_{L^{2}(\Omega \times \T^{d})} \right )
\end{aligned}
$$
 From here we obtain \eqref{propasym}$_{2}$, \eqref{propasym}$_{3}$ (along a subsequence), 
and that 
$$
g \in H^{-1} (\Omega ; L^{2} (\T^{d}))
$$

In addition, we have
$$
\begin{aligned}
\nabla_{v} g_{\eps} &= \eps \big ( \nabla_{x} + \frac{1}{\eps} \nabla_{v} \big ) g_{\eps}
- \eps\nabla_{x} g_{\eps}
\\
&= (\nabla_{x} u_{\eps}) \big (  \delta_{p}(v - \frac{x}{\eps}) - 1 \big ) 
- \eps \nabla_{x} g_{\eps}
\\
&\in_{b} L^{2}(\Omega ; M_{p} ) + H^{-2} (\Omega ; M_{p} )
\end{aligned}
$$
Therefore,
$$
\begin{aligned}
< \nabla_{v} g_{\eps} , \theta > &= \int_{\Omega} \nabla_{x} u_{\eps}(x)
\Big ( \theta(x, \frac{x}{\eps}) - \int_{\T^{d}} \theta(x,v)dv \Big ) dx
+ \eps <g_{\eps}, \nabla_{x} \theta >
\\
| <\nabla_{v} g , \theta >| &= \lim | <\nabla_{v} g_{\eps} , \theta >|
\\
&\le {\overline \lim} \left [ 
\int_{\Omega} |\nabla_{x}u_{\eps}(x) \theta (x, \frac{x}{\eps}) | dx
+
\int_{\Omega} \int_{\T^{d}} |\nabla_{x}u_{\eps}| | \theta |  dv dx
\right ]
\\
&\le C \| \theta \|_{L^{2}(\Omega \times \T^{d})}
\end{aligned}
$$
which implies
$$
\nabla_{v} g \in L^{2} (\Omega \times \T^{d})
$$
and gives using the Poincar\'e inequality the desired \eqref{propasym}$_{4}$.

To see the first property, consider a test function $\theta = \varphi \otimes \psi$
which is a tensor product of $\varphi \in C^{\infty}_{c}(\Omega)$
and  $\psi \in C^{\infty} (\T^{d})$. Then
$$
\begin{aligned}
&<\frac{1}{\eps} \big (\delta_{p} (v - \frac{x}{\eps}) - 1 \big ) , 
\varphi \otimes \psi >
\\
&\quad
=  \sum_{k \in \Z^{d}} \eps^{d} \int_{\T^{d}}\frac{1}{\eps}
\varphi (\chi_{k} + \eps w )
\Big (  \psi(w)  - \int_{\T^{d}} \psi \Big )
dw
\\
&\quad 
=
\sum_{k \in \Z^{d}} \eps^{d} \int_{\T^{d}} 
\frac{  \varphi(\chi_{k} + \eps w ) - \varphi(\chi_{k})  
- \nabla \varphi (\chi_{k}) \cdot \eps w}{\eps}
\big (  \psi(w)  - \int_{\T^{d}} \psi \big ) dw
\\
&\qquad 
+ \sum_{k \in \Z^{d}} \eps^{d} \nabla \varphi (\chi_{k}) \cdot 
\int_{\T^{d}} w \Big (  \psi(w)  - \int_{\T^{d}} \psi \Big ) dw
\\
&\quad = O(\eps) + \int_{\Omega} \nabla \varphi (x) dx \cdot
\int_{\T^{d}} w \Big (  \psi(w)  - \int_{\T^{d}} \psi \Big ) dw
\\
&\quad
\to 0
\end{aligned}
$$
as $\varphi$ is of compact support. Since
$$
\frac{1}{\eps} \big (\delta_{p} (v - \frac{x}{\eps}) - 1 \big ) 
\in_{b}  H^{-1} (\Omega ; M_{p}) = \left ( H_{0}^{1} (\Omega ; C_{p}) \right )^{*}
$$
and finite sums of tensor products $\sum_{j} \varphi_{j} \otimes \psi_{j}$ are dense
in $H_{0}^{1} (\Omega ; C_{p})$ we obtain \eqref{propasym}$_{1}$.
\end{proof}

%

\section{Homogenization of hyperbolic equations}
\label{homohyp}

In this section we consider certain homogenization problems for transport equations.
First we develop an example where the effective equation can be calculated with the help
of the double scale kinetic decomposition. Then we provide a counter-example where the 
double scale limit is not the right object to treat the effective equation.

%
\subsection{Effective equation}
\label{subseffeq}
Consider the transport equation
\begin{equation}
\label{transport}
\begin{aligned}
\dt {u_{\eps}} + a \big (\frac{x}{\eps} \big) \cdot \nabla_{x} u_{\eps} &= 0 
\\
u_{\eps} (0, x) &= U^{0}(x, \frac{x}{\eps})
\end{aligned}
\end{equation}
We assume that $a(v)$ is a $C^{1}$ vector field, periodic with period $1$,
and satisfying $\div\, a = 0$, and that the initial data 
 $U^{0} \in L^{2}(\R^{d} \times \T^{d})$ is 1-periodic in $v$
and satisfy the uniform bounds
\begin{equation}
\label{hypdata}
\int_{\R^{d}} | U^{0} (x, \frac{x}{\eps}) |^{2} dx \le C \, .
\tag {h$_{d}$}
\end{equation}
Under this hypothesis standard energy estimates for \eqref{transport} imply the 
uniform bound on solutions
\begin{equation}
\label{unifl2}
\int_{\R^{d}} | u_{\eps} (t, x) |^{2} dx 
   \le \int_{\R^{d}} | U^{0} (x, \frac{x}{\eps}) |^{2} dx 
\le C \quad \forall t > 0 \, .
\end{equation}

We introduce 
\begin{equation}
\label{dsdef}
f_{\eps} (t,x,v) = u_{\eps}(t,x) \delta_{p} \big( v -  \frac{x}{\eps}  \big)
\quad t \in \R_{+}, \, x\in \R^{d} , \, v\in\T^{d} \, ,
\end{equation}
and use Lemma \ref{lem1} to check that $f_{\eps}$ satisfies
\begin{align}
\dt {f_\ep} +a(v)\cdot\nabla_x f_\ep +\frac{1}{\ep}a(v)\cdot
\nabla_v f_\ep &=0,\quad \text{in $\cD'$}
\label{maineq}
\\
f_\ep(t=0,x,v) 
&=U^0(x,v)\;\delta_{p} (v-\frac{x}{\eps})
\label{newinit} 
\end{align}
with periodic boundary conditions ($v \in \T^d$).
The uniform bound \eqref{unifl2} implies
\begin{equation}
f_\ep \in_{b}
L^\infty \big ( [0,\ \infty) ,\ L^2(\R^d, M_{p}) \big ) \, ,
\label{bounds}
\end{equation} 
and thus by Proposition \ref{prop1}, along a subsequence if necessary,
\begin{equation}
\label{weakl}
f_{\eps} \rightharpoonup f \quad \text{weak-$\star$ in  
$L^\infty \big ( [0,\ \infty) ,\ L^2(\R^d, M_{p}) \big )$}
\end{equation}
with $f$ enjoying the improved regularity
\begin{equation}
\label{wlreg}
f \in L^\infty([0,\ \infty) ,\ L^2(\R^d \times \T^{d})) \, .
\end{equation}

Our objective is to calculate the effective limit of \eqref{transport} by computing the
hydrodynamic limit problem for the kinetic equation \eqref{maineq}-\eqref{newinit}.
Note that if $f$ satisfies a well-posed problem then this provides a complete determination
of the weak limit of $u_{\eps}$ since
$$
u_{\eps}  = \int_{\T^{d}} f_{\eps} dv \rightharpoonup \int_{\T^{d}} f dv = u
$$
We introduce 
\begin{equation}
\label{projkernel}
K=\Big\{g\in L^2(\T^d)\quad 
\Big|\ a(v)\cdot
\nabla_v g=0\ in \ \cD'\Big\}.
\end{equation}
and remark that $K$ is the space of solutions of the cell-problem obtained by the method
of multiscale asymptotic expansion \cite{BLL} for the homogenization problem 
\eqref{transport} (see \cite{HX}, \cite{E}). 
Let $P$ denote the $L^{2}$-projection operator on the kernel $K$.
We prove

\begin{theorem}
\label{hyphomthm}
Under hypothesis \eqref{hypdata} the effective limit of problem \eqref{transport}
is obtained as $u = \int_{\T^{d}} f dv$ where 
$f \in L^\infty([0,\ \infty) \, , \, L^2(\R^d \times \T^{d}))$,
$f(t,x,\cdot) \in K$ for a.e. $(t,x)$, and $f$
is the unique solution of the kinetic problem
\begin{equation}
\begin{aligned}
\dt f + (P a) \cdot \nabla_{x} f &= 0
\\
f(t=0, x,v) &= P U^{0} (x,v) \, , 
\end{aligned}
\end{equation}
where $P$ is the projection operator on the kernel $K$.
\end{theorem}

\begin{proof} Let $f_{\eps}$  and $f$ be as in \eqref{bounds}, \eqref{weakl}
and \eqref{wlreg}. The proof is split in three steps:

{\em Step 1~: The limit $f$ belongs to $K$.} 
The kernel $K$ is defined in \eqref{projkernel}.
We may consider elements of $K$ as functions of $t,\,x$ and $v$ instead of
only $v$, as $t$ and $x$ play the role of parameter in the
definition of $K$. Thus we have  
\begin{equation}
K_x=\Big\{g\in L^2(\R^d\times\T^d)\quad 
\Big|\ a(v)\cdot
\nabla_v g=0\ in\ \cD'\Big\},
\end{equation}
and
\begin{equation}
K_{t,x}=\Big\{g\in L^\infty([0,\ \infty),\ L^2(\R^d\times\T^d))\quad 
\Big|\ a(v)\cdot
\nabla_v g=0\ in\ \cD'\Big\}.
\end{equation}
We may also define all the
\[
K^p=\Big\{g\in L^p(\T^d))\quad 
\Big|\ a(v)\cdot
\nabla_v g=0\ in\ \cD'\Big\},
\]
and their extensions $K_x^p$ and $K_{t,x}^p$.

The convergence \eqref{weakl} states that for $\phi$ in 
$L^1([0,\ \infty),\ L^2(\R^d,\ C_{p}))$  we have
\[
\int_0^\infty\int_{\R^d\times\T^d} \phi(t,x,v)\,df_\ep \longrightarrow
\int_0^\infty\int_{\R^d\times\T^d} \phi(t,x,v)\,df.
\]
Take $\phi\in C^\infty_c([0,\ \infty)\times\R^d\times\T^d)$ and compute
\[
\int_0^\infty\int_{\R^d\times\T^d} a(v)\cdot\nabla_v \phi\,df_\ep=
-\ep\int_0^\infty\int_{\R^d\times\T^d}(\partial_t\phi+a(v)\cdot
\nabla_x \phi)\,df_\ep.
\]
Passing to the limit, we conclude that
\[
\int_0^\infty\int_{\R^d\times\T^d} a(v)\cdot\nabla_v \phi\,df=0.
\]
On the other hand $f \in L^\infty([0,\ \infty), L^2(\R^d\times\T^d))$, so 
$f\in K_{t,x}$.

\medskip
{\em Step 2~: The limit equation.} 
Consider a function $\phi \in K_x$. We wish to mollify $\phi$ and use it
as a test function in the weak form of \eqref{maineq}. Since $K_x$ depends
only parametrically in $x$, we may select $\phi$ to be compactly
supported in $x$.

Take $H(x) \in C^\infty_c(\R^d) $ with $\int_{\R^d} H(x)\,dx=1$, and
 $\bar H(v) \in C^\infty_c ((0,\ 1)^d)$ with
$\int_{\T^d} \bar H(v)\,dv=1$. 
For any $\phi$ define
\[\begin{split}
&\phi_{n}=\int_{\R^d} n^d\,H\left(n(x-y)\right)\;\phi(y,v)\,dy,\\
&\phi_{n,m}=\int_{\T^d}
\bar H_m(v-\eta)\;\phi_n(x,\eta)\,d\eta,\\
\end{split}
\]
with $\bar H_m(v)=m^d\sum_{k\in \Z^d} \bar H(m(v+k))$, periodic
and well
defined for all $m$ as $\bar H(m v)$ is compactly supported in
$(0,\ 1/m)^d$.

Then for any $\phi\in K_x$  we have
$$
\int_{\R^{d}\times \T^{d}} \phi(y,\eta) a(\eta) \cdot
\nabla \bar H_{m} (v-\eta) H_{n} (x-y) \, dy d\eta =0 \, .
$$
Thus, for $a$ Lipshitz continuous,
$$
\begin{aligned}
\int_{\R^d\times\T^d} &a(v)\cdot\nabla_v\phi_{n,m}\,df_\ep
\\
&=
\int_{\R^d\times\T^d}\!\!\int_{\T^d} (a(v)-a(\eta))\cdot \nabla \bar H_m
\left(v-\eta\right) \phi_n(x,\eta)\,d\eta\,df_\ep
\\
&=  \int_{\R^d\times\T^d}\!\!\int_{\T^d}\int_0^1 
\big(\zeta\cdot\nabla a(v-(1-t)\zeta)\big)\,
\\
&\qquad\qquad\qquad\qquad
\cdot \sum_{k\in\Z^{d}} m^{d+1}\nabla \bar H\left(m(\zeta +k)\right)
\phi_n(x,v-\zeta)\,dt\,d\zeta\,df_\ep.\\
\end{aligned}
$$
Notice that $m^{d+1}\zeta\otimes \nabla\bar H(m\zeta)$ converges in 
the sense of distributions toward $C\,(Id)\,\delta$ with $C$ a numerical constant.
Moreover thanks to \eqref{bounds}, and to the fact that $\phi_n\in 
L^2(\T^d,\ C_c(\R^d))$ and $\nabla a\in C (\T^d)$, we may pass to
the limit in $m$ in the previous equality and find
\[\begin{split}
\lim_{m\rightarrow\infty}
\int_{\R^d\times\T^d} &a(v)\cdot\nabla_v\phi_{n,m}\,df_\ep=C\,
\int_{\R^d\times\T^d} {\rm div}\, 
a(v)\,\phi_n(x,v)\,df_\ep=0,
\end{split}
\]
as $a$ is divergence free. 
Multiplying \eqref{maineq} by $\phi_{n,m}$ and taking first
the limit $m\rightarrow \infty$ and then the limit $\ep\rightarrow 0$, 
we find that for any $\phi\in K_x$ compactly supported in $x$ we have
\[
\partial_t \int_{\R^d\times \T^d} \phi_n f\,dx\,dv -\int_{\R^d\times\T^d}
a(v)\cdot\nabla_x\phi_n\,f\,dx\,dv=0.
\]
This relation can now easily be extended by approximation to any $\phi \in K_x$.

Let us denote $\bar a$ the orthogonal projection on $K$ of $a$. The new 
function $\bar a$ belongs to $L^\infty(\T^d)$ as the projection
operator $P$ is continuous on every $L^p(\T^d)$ for all $1\leq p\leq\infty$, 
but does not necessarily
have any further
regularity, Lipschitz for instance (see the appendix where we recall
the basic properties of $P$). Now as $f\in K_{t,x}$ and
$\nabla_x\phi_n\in K_x$,
then $\nabla_x\phi_n \,f\in K^1_{t,x}$ and consequently
\[
\frac{d}{dt} \int_{\R^d\times \T^d} \phi_n f\,dx\,dv -\int_{\R^d\times\T^d}
\bar a(v)\cdot\nabla_x\phi_n\,f\,dx\,dv=0.
\] 
On the other hand, the projection operator $P$ may be trivially extended on
$K_{t,x}$ from $K$ as $t$ and $x$ are only parameters and 
of course it commutes 
with derivatives in $t$ or $x$.  Now, for any $\phi\in L^2(\R^d\times\T^d)$
\[
\int_{\R^d\times\T^d} (\phi-P\phi)\,f\,dx\,dv=0,
\int_{\R^d\times\T^d} \bar a\cdot\nabla(\phi_n-P\phi_n)\,f\,dx\,dv=0.
\]
Finally for any $\phi\in L^2(\R^d\times\T^d)$, we have that
\[
\partial_t \int_{\R^d\times \T^d} \phi_n f\,dx\,dv -\int_{\R^d\times\T^d}
\bar a(v)\cdot\nabla_x\phi_n\,f\,dx\,dv=0.
\] 
This implies that $f \in K_{t,x}$ is a solution in the sense of distribution to 
\begin{equation}
\partial_t f+\bar a(v)\cdot\nabla_x f=0.\label{limiteq}
\end{equation}

\medskip
{\em Step 3~: Conclusion.} Let us begin with the identification of the
initial value $f(t=0)$ which has a sense since $\partial_t f\in L^\infty([0,\
\infty),\ H^{-1}(\R^d,\ L^2(\T^d)))$ because of \eqref{limiteq} and as
$f\in L^\infty([0,\
\infty),\ L^2(\R^d\times \T^d)))$. For every $\phi\in K_x$, as
\begin{equation}
\frac{d}{dt} \int_{\R^d\times \T^d} \phi_n df_\ep =\int_{\R^d\times\T^d}
a(v)\cdot\nabla_x\phi_n\,df_\ep,
\label{inter}
\end{equation}
then $\int_{\R^d\times \T^d} \phi_n df_\ep(t,.,.)$ has a limit as 
$t\rightarrow 0$ and this limit is, thanks to \eqref{newinit}
\[
\int_{\R^d} \phi_n(x,x/\ep)\,U^0(x,x/\ep)\,dx. 
\]
Moreover because of \eqref{inter}, we may pass to the limit in $\ep$
and  deduce that
\[
\int_{\R^d\times \T^d} \phi_n f(t,x,v)\,dxdv
\mathop{ \; \; \longrightarrow \; \;}_{t\to 0}
\int_{\R^d\times \T^d} \phi_n f(0,x,v)=
\int_{\R^d\times \T^d} \phi_n U^0(x,v).
\]
On the other hand we of course have for any $\phi\in L^2$ as $f\in K_{t,x}$
$$
0=\int_{\R^d\times \T^d} (\phi_n-P\phi_n) f(t,x,v)\,dxdv
\mathop{ \; \; \longrightarrow \; \;}_{t\to 0}
\int_{\R^d\times \T^d} (\phi_n-P\phi_n) f(0,x,v).
$$
Combining the last two equalities we get that
\begin{equation}
f(t=0,x,v)=P\,U^0(x,v).
\label{initf}
\end{equation}

Finally, we notice that Eq. \eqref{limiteq} combined with \eqref{initf} has
a unique solution in the space of distribution, through standard arguments of
kinetic theory and as, even though $\bar a$ is only bounded, it does not 
depend on $x$. Therefore any extracted subsequence of
 $f_\ep$ has only one possible limit and the whole sequence $f_\ep$ converges
toward the solution of \eqref{limiteq} with \eqref{initf}.
\end{proof}

\bigskip
\noindent
{\sc Examples.}
We calculate the equation for the double scale limit $f$ and the associated 
effective equation for certain examples, always within the framework of
\eqref{transport}.

1. First consider the case that $a(v)$ is ergodic. 
Then 
$$
\begin{aligned}
K &= \{ g \in L^2(\T^d) : g = const. \}
\\
P_K g &= \int_{\T^d} g dv = : \overline{ g }
\end{aligned}
$$
The equation for $f$ becomes
$$
\del_t f + \overline{ a} \cdot \nabla_x f = 0
$$
and of course $u = \int_{\T^d} f dv$.

2. Consider next the homogenization problem
$$
\begin{aligned}
\del_t u_\eps + b \big ( \frac{x_2}{\eps} \big ) \del_{x_1} u_\eps &= 0
\\
u_\eps (0, x_1 , x_2) &= U^0 (x_1, x_2, \frac{x_1}{\eps}, \frac{x_2}{\eps}) 
\end{aligned}
$$
where $u_\eps = u_\eps (t,x)$, $x = (x_1, x_2) \in \R^2$, 
and the vector field $a(x_1, x_2) = (b(x_2), 0)$ corresponds
to a shear flow with $b(x_2)\neq 0$ for $a.e.\ x_2$. We compute
$$
\begin{aligned}
K &= \{ g \in L^2(\T^2) : b(v_2) \del_{v_1} g = 0 \}
= \{ g = \psi (v_2) \; \big | \;  \forall \psi \in L^2(\T^1) \}
\\
P_K g &= \int_0^1 g(v_1, v_2) dv_2 
\end{aligned}
$$
Since $f \in K$ we conclude that $f = f(t, x_1, x_2, v_2)$ and
satisfies the problem
$$
\begin{aligned}
\del_t f + b(v_2) \del_{x_1} f &= 0
\\
f (0, x_1, x_2 , v_2) &= P_K U^0 = \int_0^1 U^0 (x_1, x_2, v_1, v_2) dv_1
\end{aligned}
$$
The weak limit $u = \int_{\T^2} f $  satisfies the 
integrated equation.

3. It is possible to give a more general framework for the situation
   of the previous example.
Suppose that the divergence free vector field $a$ is such that the
   following description of $K$ is true: There exist
   functions $\xi_1,\ldots \xi_N$, $N\leq d$ from $\T^d$ to
   $\R$. These functions are
   local coordinates in the sense that they may be
   completed  by $\xi_{N+1},\ldots \xi_d$ and that the change of
   coordinates $v$ to $(\xi_1(v),\ldots, \xi_d(v))$ is a $C^1$
   diffeomorphism from $\T^d$ to some domain $O\subset \R^d$. And finally
\[
K=\{\psi(\xi_1(v),\ldots,\xi_N(v))\;|\ \forall \psi\in L^2(O) \}.
\]
For instance in dimension $d=2$, as ${\rm div}\;a=0$, there
is always $\xi:\ \R^d\rightarrow \R$ such that $a=\nabla^\perp
\xi$. Now if in addition
 $\xi$ is a periodic regular function with
$\nabla \xi(v)\neq 0$ for all $v$, which is a non trivial assumption,
then $K$ is exactly the set of functions $\psi(\xi)$. 

In that case, we may define $g(t,x,\xi)=f(t,x,V(\xi_1,\ldots,\xi_d))$ 
with $V$ the inverse change of variables. Then $g$ does not depend on
$\xi_{N+1},\ldots , \xi_d$ and it simply satisfies
\[
\partial_t g+b(\xi_1,\ldots,\xi_N)\cdot\nabla_x g=0,
\]
with $b(\xi)=\bar a(V(\xi))=\int a(V(\xi))d\xi_{N+1}\ldots d\xi_d$.

4. Notice now that the kernel $K$, endowed with the usual $L^2$ scalar
product, is a Hilbert space and so
that the kernel $K$ admits an orhonormal basis $\{ \psi_k (v) \}$,
possibly countable. Since $f \in K$ it will be given in a Fourier
expansion
$$
f = \sum_{k=1}^\infty m_k (t,x) \psi_k (v) 
\quad \text{where $m_k = < f , \psi_k >$}.
$$
Moreover, we see that
$$
<P_K (a) f , \psi_k > = < P_K (a f) , \psi_k > = < af , \psi_k >
$$
and one computes that the set of moments $m_k$ satisfies the initial 
value problem
$$
\begin{aligned}
&\del_t m_k + \sum_{j=1}^d  \Big ( \sum_{n=1}^\infty
<a_j \psi_n , \psi_k > \frac{ \del m_n}{\del x_j} \Big ) = 0
\\
&m_k (0, x) = <P_K U^0 (x, \cdot ) , \psi_k > 
= \int_{\T^d} U^0(x,v) 
\overline{ \psi_k (v) } dv
\end{aligned}
$$
As the wave speed $a$ is real, 
$<a_j \psi_n , \psi_k > = <\psi_n , a_j \psi_k >$, and the system of moments 
is an infinite symmetric hyperbolic system.

%

\subsection{The multiscale case: A counterexample}
\label{secce}
A natural extension of the previous analysis is to deal with transport
coefficients depending on more than one scale. Consider for example
the equation
\begin{equation}
\partial_t u_\eps+a_\eps\cdot\nabla_x u_\eps=0,
\label{multiscaletransport}
\end{equation}
with $a_\eps=a(x,x/\eps)$ and $a(x,v)$ a Lipschitz function,
or even with $a_\eps=a(x,x/\eps,x/\eps^2)$ (or with as many
scales as one cares to introduce). Assume again that
$\hbox{div}_v a(x,v)=0$ and $\hbox{div}_w a(x,v,w)=0$.

Is it possible to derive an equation for the double scale limit (or
for the triple scale limit when $a(x,x/\eps,x/\eps^2)$) in the case of
\eqref{multiscaletransport}?

In fact, it is relatively easy to show that the previous approach does not
work! Everything goes as before in the beginning; upon defining
\[
f_\eps(t,x,v,w)=u_\eps(t,x)\,\delp(v-x/\eps)\,\delp(w-v/\eps),
\]
as in paragraph \ref{ss13}, one simply obtains the generalized kinetic
equation
\[
\partial_t
f_\eps+\nabla_x\cdot(a(x,v,w)f_\eps)+\frac{1}{\eps}a\cdot\nabla_v
f_\eps+\frac{1}{\eps^2}a\cdot\nabla_w f_\eps=0.
\]

However it is not always possible to derive a well posed problem for the
hydrodynamic limit, even in the simple setting where $a$ depends only on
$x,v$ and the equation for $f_\eps(t,x,v)$ is
\[
\partial_t
f_\eps+\nabla_x\cdot(a(x,v)f_\eps)+\frac{1}{\eps}a(x,v)\cdot\nabla_v
f_\eps=0.
\]
Indeed the only information that we have is that any limit $f$ belongs
to the kernel which now depends on $v$ and $x$
\[
K=\{f\in L^2(\R^d\times\T^d)\;|\; a(x,v)\cdot\nabla_v f=0\}.
\]
On the other hand, when projecting the equation on $K$, it is not possible to
handle the term with the $x$ derivative as projection on $K$ and
differentiation in $x$ no longer commute. This is associated to the
possibility that the dimensionality of $K$ may vary with $x$.

In addition, at the level of the double scale limit, this is not a mere
technical problem, rather the double scale limit is in general not unique and
depends on the choice of the extracted subsequence in $\eps$.

This can be simply seen for the problem
\begin{equation}
\begin{split}
&\partial_t u_\eps +a(x)\cdot\nabla_x u_\eps=0,\quad t\in\R_+,\
x\in\R^2,\\
&u_\eps(t=0,x)=U^0(x,x/\eps).
\end{split}\label{c-ex}
\end{equation}
The oscillations are due only to the initial data as the transport
coefficient no longer depends on $\eps$, and $u$ the weak limit
of $u_\eps$ satisfies the same equation. Take now
\[
a_1(x)=1,\quad a_2(x)=x_1\,\charf_{0\leq x_1\leq 1}+1\,\charf_{x_1>1},
\]
so that $a$ is Lipschitz and divergence free, and select the initial data
\[
U^0(x,v)=K(x)\,L(v_2),
\]
with $K$ and $L$ two $C^\infty$ functions, $L$ periodic on $\R$ of
period $1$ and with zero average, and $K$ compactly supported with support in $x_1$ in $[-1,\
-1/2]$.

As the average of $U^0(x,v)$ in $v$ vanishes for all $x$, the weak
limit $u$ of $u_\eps$ is uniformly $0$. 

For $1\leq t\leq 3/2$, the support in $x_1$ of the solution $u_\eps$ 
is entirely in the interval $[0,\ 1]$. Therefore any double scale
limit $f$ should satisfy
\[
\partial_{v_1} f+x_1\partial_{v_2}f=0.
\]
It is easy to check that the only $L^2$ solutions to this last equations are the functions
which depend only on $x$ and not on $v$. Therefore for $1\leq t\leq
3/2$ the double scale limit is equal to $u$, {\rm i.e.} uniformly
vanishes.

Let us finally compute the double scale limit for $t>2$ and check that
it does not vanish. For that introduce the characteristics $X(t,x)$
\[
\partial_t X(t,x)=a(X(t,x)),\quad X(0,x)=x.
\]
We need the
characteristics only for those $x$ which belong to the support of
$U^0$, that is for $-1<x_1<-1/2$.

As $a_1=1$, we simply have
\[
X_1(t,x)=x_1+t.
\]
As to $X_2$, as long as $X_1<0$ or $t<t_0=-x_1$ (remember $x_1\in[-1,\
-1/2]$) it is equal to $x_2$.
For $t_0<t<t_1=1-x_1$ (corresponding to $X_1$ in $[0,\ 1]$), we have
\[
\partial_t X_2=X_1=x_1+t.
\] 
As a consequence 
\[
X_2(t_1)=x_2+x_1\,(t_1-t_0)+\frac{t_1^2}{2}-\frac{t_0^2}{2}=x_2+x_1+\frac{(1-x_1)^2}{2}-
\frac{x_1^2}{2}=x_2+1/2.
\]
After $t>t_1$, $\partial_t X_2=1$ and so
\[
X_2(t)=x_2+1/2+t-t_1=x_2+1/2+t-1+x_1=x_2+X_1-1/2.
\] 
With this, the solution $u_\eps$ is given for $t>2$ by
\[\begin{split}
u_\eps(t,x)&=u_\eps(0,x_1-t,x_2-x_1+1/2)\\
&=K(x_1-t,x_2-x_1+1/2)\, L(x_2/\eps-x_1/\eps+1/2\eps).
\end{split}\]
For every $\alpha\in [0,\ 1]$, choose a subsequence $\eps_n$ such that
$1/2\eps_n$ converges to $\alpha$ modulo $1$. Then the double scale
limit associated to this subsequence is the function
\[
f(t,x,v)=K(x_1-t,x_2-x_1+1/2)\, L(v_2-v_1+\alpha).
\]
Instead of one unique limit, we obtain a whole family which clearly
indicates the ill-posedness of the problem at the level of the double
scale limit.

%

\section{Transport of oscillating fine-scale}
\label{sectr}

An interesting question that can be studied using the techniques 
developped in section \ref{homohyp} is the problem of transport of an oscillatory fine-scale structure 
under a divergence-free vector field. Consider the homogenization problem
\begin{equation}
\label{tos}
\begin{aligned}
\del_{t} u_{\eps} + \nabla_{x} \cdot A(t, x, \frac{\varphi(t,x)}{\eps}) u_{\eps} &= 0
\\
u_{\eps} (0, x) = U^{0} (x, \frac{x}{\eps})
\end{aligned}
\end{equation}
where $\varphi : \R^{d} \times \R \to \R^{d}$ is a $C^{2}$ map describing the fine scale of
oscillations  that satisfies for some $c >0$
\begin{equation}
\label{hypso}
\begin{aligned}
\text{$\varphi (\cdot, t)$ is surjective and invertible for $t$ fixed}
\\
\det \nabla (\varphi^{-1} (\cdot , t)) \ge c > 0
\\
\varphi (\cdot , 0 ) = id
\end{aligned}
\tag {h$_{so}$}
\end{equation}
and $a_{\eps} = A(t, x, \frac{\varphi(t,x)}{\eps})$ is a divergence free field.

The latter is guaranteed provided $A(t,x,v)$ is a $C^{1}$ vector field $1$-periodic in $v$
such that
\begin{equation}
\label{hyptvf}
\begin{aligned}
\nabla_{x} \cdot A(t,x,v) &= 0
\\
\tr \big ( \nabla_{v} A \nabla_{x} \varphi \big ) 
= \sum_{i,j = 1}^{d} \frac{\del A_{i}}{\del v_{j}} \frac{\del \varphi_{j}}{\del x_{i}}
&=0
\end{aligned}
\tag {h$_{tvf}$}
\end{equation}
The initial data 
 $U^{0} \in L^{2}(\R^{d} \times \T^{d})$ are 1-periodic in $v$
and satisfy the uniform bounds
\begin{equation}
\label{todata}
\int_{\R^{d}} | U^{0} (x, \frac{x}{\eps}) |^{2} dx \le C \, .
\tag {h$_{d}$}
\end{equation}
Then standard energy estimates imply that solutions of \eqref{tos}
satisfy the uniform bound
\begin{equation}
\label{uniforml2}
\int_{\R^{d}} | u_{\eps} (t, x) |^{2} dx 
   \le \int_{\R^{d}} | U^{0} (x, \frac{x}{\eps}) |^{2} dx 
\le C  \, .
\end{equation}
Our objective is to calculate an effective equation for the weak limits of $\{u_{\eps}\}$.
The counterexample of section \ref{secce} indicates that we can not expect to do that in full
generality. A more precise statement of what will be achieved is that we will identify
conditions on the vector field $A$
and the structure function $\varphi$ under which an effective equation is calculated.

\subsection{Reformulation via a kinetic problem}
We introduce the "kinetic function" 
\begin{equation}
\label{kinf}
f_{\eps} = u_{\eps} (t,x) \delta_{p} \big (v - \frac{\varphi (t,x) }{\eps} \big ) \, ,
\end{equation}
which is well defined (see section \ref{ss12}) as a measure. Due to the identities 
in lemma \ref{lem2} of section \ref{ss12}, it is possible to transform
the homogenization problem \eqref{tos} into a hyperbolic limit for a kinetic initial value problem:

\begin{lemma}
If $u_{\eps}$ a weak solution of \eqref{tos} then
$f_{\eps}$ in \eqref{kinf} verifies in $\cD'$ the  kinetic problem
\begin{equation}
\label{kintos}
\begin{aligned}
\del_{t} f_{\eps} + \nabla_{x} \cdot (A f_{\eps}) + \nabla_{v} ( \frac{1}{\eps} B f_{\eps})
&= 0 \, ,
\\
f_{\eps}(0,x,v) &= U^{0} (x,v) \delta_{p} (v - \frac{x}{\eps}) \, ,
 \end{aligned}
 \end{equation}
 where the vector field $B(t,x,v)$, defined by
 \begin{equation}
\label{trvf}
B_{i}  =  \big (\del_{t} + A \cdot \nabla_{x} \big ) \varphi_{i}
\, , \quad i = 1 , ... , d \, ,
\end{equation}
is $1$-periodic in  $v$ and divergence-free,
$\nabla_{v} \cdot B = 0$.
\end{lemma}

\begin{proof} 
Let $u_{\eps}$ be a weak solution of \eqref{tos}. By \eqref{hyptvf},  $B$ defined in \eqref{trvf}
is $1$-periodic in $v$ and divergence free.
We consider a test function
$\theta \in C_{c}^{1} \big ( [0, \infty ) \times \R^{d} ; C^{1} (\T^{d}) \big)$
and compute
$$
\begin{aligned}
\int_{0}^{\infty} \int_{\R^{d}} \int_{\T^{d}}
&\big [   \del_{t} \theta    + \nabla_{x} (A \cdot \theta )
+ \frac{1}{\eps} \nabla_{v} ( B \cdot\theta) \big ] df_{\eps}
\\
&=  \int_{0}^{\infty} \int_{\R^{d}}
 u_{\eps}(x,t) \Big (   \del_{t} \theta    + A \cdot \nabla_{x} \theta 
+ \frac{1}{\eps} B \cdot \nabla_{v} \theta \Big )  (t,x,\frac{\varphi}{\eps})
dxdt
\\
&= \int_{0}^{\infty} \int_{\R^{d}}
u_{\eps} \left [   \del_{t} \big (\theta (t,x,\frac{\varphi}{\eps})  \big )  
+ A (t,x,\frac{\varphi}{\eps}) \cdot \nabla_{x} \big ( \theta (t,x,\frac{\varphi}{\eps})  \big )  
\right ] 
dxdt
\\
&= - \int_{\R^{d}} U^{0} (x, \frac{x}{\eps}) \theta(0,x,\frac{x}{\eps}) dx
\\
&= - \int_{\R^{d}} \int_{\T^{d}}  \theta(0,x,v) df_{\eps} (0,x,v) \, ;
\end{aligned}
$$
that is, $f_{\eps}$ is a weak solution of \eqref{kintos}.
\end{proof}

An alternative, albeit formal, derivation of \eqref{kintos} may be obtained by studying
characteristics. The characteristic curve of \eqref{tos} emanating from the point $y$
is defined by 
$$
\begin{cases}
\frac{dx}{dt} = A(x,t, \frac{\varphi}{\eps}) & \\
x(0,y) = y & \\
\end{cases}
$$
and is denoted by $x = X (t ; y)$. 
Along such curves we have
$$
\frac{d}{dt} \Big (\frac{\varphi}{\eps} \Big ) = \frac{1}{\eps} 
\Big (   \varphi_{t} + A(x,t,\frac{\varphi}{\eps})\cdot \nabla_{x} \varphi \Big )
$$
The two equations together can be embedded into the system of ordinary
differential equations
\begin{equation}
\label{charsys}
\begin{cases}
\frac{dx}{dt} = A(t,x,v) & \\
\frac{dv}{dt} = \frac{1}{\eps} 
\Big ( \varphi_{t} + A(x,t,v)\cdot \nabla_{x} \varphi\Big )
\end{cases}
\end{equation}
in the following sense: If $(Y(t; y,u), U(t; y, u))$ is the solution of \eqref{charsys}
emanating from the point $(y,u)$ then
$$
\begin{aligned}
X(t ; y) &= Y ( t ; y , \varphi (y,0)) 
\\
\frac{\varphi( X(t;y), t)}{\eps} &= U (t ; y , \varphi (y,0))
\end{aligned}
$$
Note that \eqref{kintos}$_{1}$ is precisely the Liouville equation associated to the characteristic
system \eqref{charsys}.

\subsection{Conditions leading to an effective equation}
Our next goal is to derive an effective equation for the hydrodynamic limit
of \eqref{kintos}. We first show that under hypothesis \eqref{hypso} 
the definition \eqref{kinf} still induces
good properties for the weak limit points of $\{ f_{\eps} \}$.

\begin{lemma}
\label{lem42}
Under hypotheses \eqref{hypso} and \eqref{uniforml2},
\begin{align}
\label{pro1}
\delta_{p} \big ( v - \frac{\varphi(x,t))}{\eps} \big ) \rightharpoonup 1 
\qquad \text{in $\cD'$}
\\
\label{pro2}
f_{\eps} \in_{b}  L^{\infty} \big ( (0,\infty) ; L^{2} (\R^{d}  ; M_{p}) \big )
\end{align}
and, along a subsequence (if necessary),
\begin{align}
\label{pro3}
f_{\eps} &\rightharpoonup f 
\qquad \text{weak-$\star$ in 
$ L^{\infty} \big ( (0,\infty) ; L^{2} (\R^{d}  ; M_{p}) \big )$}
\\
\label{pro4}
&\text{with} \; 
f \in  L^{\infty} \big ( (0,\infty) ; L^{2} (\R^{d}  \times \T^{d}) \big )
\end{align}
\end{lemma}

\begin{proof}
For $\theta \in C_{c}^{\infty} \big (   (0,\infty) \times \R^{d} ; C^{\infty} (\T^{d}) \big )$
we have
$$
\begin{aligned}
< \delta_{p} \big ( v - \frac{\varphi(x,t))}{\eps} \big ) , \theta >
&=
\int_{0}^{\infty} \int_{\R^{d}} \theta (x,t, \frac{\varphi(x,t))}{\eps}) dx dt
\\
&=
\int_{0}^{\infty} \int_{\R^{d}} \theta ( \varphi^{-1}(y,t), t, \frac{y}{\eps}) 
|\det \nabla_{y}(\varphi^{-1})| dy dt
\\
&\to
\int_{0}^{\infty} \int_{\R^{d}}\int_{\T^{d}} \theta(t,x,v) dv dx dt
\end{aligned}
$$
Note next that for $(x,t)$ fixed 
$$
\| f_{\eps} (x,t,\cdot)\|_{M_{p}} = |u_{\eps} (x,t)|
$$
and thus \eqref{pro2} and \eqref{pro3} follow from \eqref{uniforml2}. Finally,
$$
\begin{aligned}
|<f_{\eps}, \theta>|  &\le 
\| u_{\eps}\|_{L^{\infty}(L^{2})} \, 
\int_{0}^{\infty} \left (
\int_{\R^{d}} \big | \theta (x,t, \frac{\varphi(x,t))}{\eps}) \big |^{2} dx 
\right)^{\frac{1}{2}} dt
\\
&\le
C \int_{0}^{\infty} \left (
\int_{\R^{d}} \big | \theta (x,t, \frac{\varphi(x,t))}{\eps}) \big |^{2} dx 
\right)^{\frac{1}{2}} dt
\\
&\to C \int_{0}^{\infty} \left ( 
\int_{\R^{d}} \int_{\T^{d}} \big | \theta (x,t, v) \big |^{2} dx dv
\right)^{\frac{1}{2}} dt
\end{aligned}
$$
and \eqref{pro4} follows.
\end{proof}

\begin{remark}
A hypothesis of the type of \eqref{hypso} is essential for the validity 
of \eqref{pro1} and accordingly for \eqref{pro4}. For instance, in the extreme case that $\varphi$
is a constant map, $\varphi(x,t) \equiv c$, it is possible by choosing appropriate sequences
$\eps_{n} \to 0$ to achieve any weak limit
$$
\delta_{p} \big ( v - \frac{c}{\eps_{n}}) \rightharpoonup  \delta_{p} (v - v_{o})
\quad \text{with any $0 < v_{o} < 1$}.
$$
The regularity of $f$ is then no better than the regularity of $\{ f_{\eps} \}$ and
\eqref{pro4} is of course violated.
\end{remark}

We conclude by providing a formal derivation of an effective equation. Consider the
$\eps \to 0$ limit of \eqref{kintos}-\eqref{trvf} and recall that, by \eqref{uniforml2}
and lemma \ref{lem42}, we have $f_{\eps} \rightharpoonup f$ as in \eqref{pro3} and \eqref{pro4}.
Define the set
$$
K_{t,x} = \big \{  g \in L^{\infty} ( (0, \infty) ; L^{2}(\R^{d} \times \T^{d})) \; \Big | \;
               B(t,x,v) \cdot \nabla_{v} g = 0  \big \} \, .
$$
The set $K = N(\cB)$ is the null space of the operator $\cB := B\cdot \nabla_{v}$ and in general it will 
depend on $(t,x)$. We will derive the effective equation under the hypothesis
\begin{equation}
\label{hypmain}
N(\cB) \; \text{ is independent of $(x,t)$}
\tag {H}
\end{equation}
Then we have the decomposition
$$
L^{2} = N(\cB) \oplus  \overline{R(\cB^{T})} = K \oplus K^{\perp}
$$
and the spaces remain the same for any point $(x,t)$. Let $P$ denote the $L^{2}$-projection on the
set $K$. Any $\theta \in L^{2}(\T^{d})$ can be decomposed as
$$
\theta = P\theta + (I-P)\theta =: \psi + \phi
$$
Moreover the differentiation operators $\del_{t}$ and $\nabla_{x}$ commute with the 
projector $P$.

For $\psi \in N(\cB)$, using $\nabla_{x} \cdot A = \nabla_{v} \cdot B = 0$,
we derive from \eqref{kintos} that 
\begin{equation}
\label{interm}
\del_{t} < f , \psi > + \nabla_{x} <Af , \psi> = 0
\end{equation}
where the brackets denote the usual inner product in $L^{2}(\T^{d})$. One
easily sees that
$f(t,x, \cdot) \in K$ for a.e. $(t,x)$. Given $\theta \in L^{2}(\T^{d})$ 
let $\psi = P \theta$. Then
$$
<f , \theta > = <f , \psi >
$$
and
$$
<A_{i} f, \psi> = < P(A_{i} f) , \psi> = <P(A_{i} f) , \theta>
$$
Since $f \in K$ we  have $P(A_{i} f) = P(A_{i}) f$ and we conclude that \eqref{interm}
can be expressed in the form
$$
\del_{t} < f , \theta > + \nabla_{x} \cdot <P(A) f , \theta > = 0  \, , 
\quad \theta \in L^{2}(\T^{d}) \, .
$$
The effective equation thus takes the form
\begin{equation}
\label{effto}
\del_{t} f + \nabla_{x} \cdot ( PA) f = 0 \, .
\end{equation}

The above derivation of \eqref{effto} is formal and is based on hypothesis \eqref{hypmain},
which is quite restrictive especially when viewed together with \eqref{hyptvf} that has to be
satisfied simultaneously. We view this equation as a theoretical framework of when an effective
equation can be computed. To derive it rigorously one needs an analysis as in Theorem \ref{hyphomthm}
and we will not pursue the details here. The counterexample in section \ref{secce} indicates
that the hypothesis \eqref{hypmain} is essential.

We list two examples that can be viewed under the above framework. First,
the homogenization problem \eqref{transport} is a special case of \eqref{tos}
with the obvious identifications. A second example is given by the problem
\begin{equation}
\label{special}
\begin{aligned}
\del_{t} u_{\eps} + a(x) \cdot \nabla_{x} u_{\eps} &= 0
\\
u_{\eps} (x, 0) &= U (x, \frac{x}{\eps})
\end{aligned}
\end{equation}
where $a$ is a divergence free field, $\nabla_{x} \cdot a = 0$.
Define $\varphi(t,x)$ to be the backward characteristic emanating 
from the point $x$.
Then $\varphi = (\varphi_{1}, ... , \varphi_{n})$ satisfies 
$$
\begin{aligned}
\del_{t} \varphi_{i} + a(x) \cdot \nabla_{x} \varphi_{i} &= 0
\\
\varphi_{i} (0,x) = x_{i}
\end{aligned}
$$
The problem \eqref{special} fits under the framework of \eqref{tos} under
the selections
$$
A(t,x, v) = a(x) \, . \quad
B_{i} (t,x,v) = \big ( \del_{t} + a(x) \cdot \nabla_{x} \big ) \varphi_{i} = 0
\, .
$$
The kinetic equation for 
$f_{\eps} = u_{\eps} \delta_{p} \big ( v - \frac{\varphi}{\eps}\big)$
becomes 
$$
\begin{aligned}
\del_{t} f_{\eps} + a(x) \cdot \nabla_{x} f_{\eps } &=0
\\
f_{\eps} (0,x,v) = U (x,v) \delta_{p} \big ( v - \frac{x}{\eps}\big)
\end{aligned}
$$
while the limiting $f$ satisfies the same transport equation with initial condition
$f(0,x,v) = U(x,v)$. Hence, it is computed explicitly by
$$
f(t,x,v) = U (\varphi(t,x),v) \, .
$$

%

\section{Enhanced diffusion}
\label{endiff}

In this section we study the enhanced diffusion problem
\begin{equation}
\label{enhdiff}
\begin{aligned}
\del_{t} u_{\eps} + \frac{1}{\eps} a(x, \frac{x}{\eps}) \cdot \nabla_{x} u_{\eps}
&= \alpha \triangle_{x} u_{\eps}
\\
u_{\eps} (0, x) = U^{0} (x, \frac{x}{\eps})
\end{aligned}
\end{equation}
where $a(x,v)$ is a Lipshitz vector field periodic (with period 1) in $v$ that satisfies
\begin{equation}
\label{vecf}
\nabla_{x}\cdot a = \nabla_{v} \cdot a =0 \, , \quad \int_{\T^{d}} a(x,v) dv = 0 \, ,
\tag {h$_{vf}$}
\end{equation}
 $\alpha > 0$ is constant and $U^{0} \in L^{2} (\R^{d} \times \T^{d})$.
We use this as an example to develop the methodology of section \ref{msde}.
For previous work and a commentary on the significance of this problem 
we refer to Avellaneda-Majda \cite{AM},
Fannjiang-Papanicolaou \cite{FP} and references therein.
It is assumed that the initial data oscillates 
at the scale $\eps$ and satisfy the uniform bound
\begin{equation}
\label{diffdata}
\int_{\R^{d}} | U^{0} (x, \frac{x}{\eps}) |^{2} dx \le C \, .
\tag {h$_{d}$}
\end{equation}
Standard energy estimates then imply the 
uniform bounds
\begin{equation}
\label{unifdiff}
\int_{\R^{d}} | u_{\eps} (t, x) |^{2} dx 
+ \alpha \int_{0}^{t} \int_{\R^{d}} |\nabla u_{\eps}|^{2} dxdt
   \le \int_{\R^{d}} | U^{0} (x, \frac{x}{\eps}) |^{2} dx 
\le C 
\end{equation}
for solutions of \eqref{enhdiff}.

We introduce the kinetic decomposition
\begin{equation}
\label{diffdsdef}
f_{\eps} (t,x,v) = u_{\eps}(t,x) \delta_{p} \big( v -  \frac{x}{\eps}  \big)
\quad t \in \R_{+}, \, x\in \R^{d} , \, v\in\T^{d} \, ,
\end{equation}
and use Lemma \ref{lem1} to obtain that $f_{\eps}$ satisfies the transport-diffusion equation
\begin{equation}
\label{maindiff}
\begin{aligned}
\dt {f_\ep} + \frac{1}{\eps} a(x, v)\cdot\nabla_x f_\ep 
&+ \frac{1}{\ep^{2}} 
\big ( a(x,v) \cdot \nabla_{v} f_{\eps} - \alpha \triangle_{v} f_{\eps} \big )
\\
 &= \alpha \triangle_{x} f_{\eps} 
 + \frac{2\alpha}{\eps} \nabla_{x} \cdot \nabla_{v} f_{\eps}
\, , \qquad \text{in $\cD'$}
\\
f_\ep(t=0,x,v) 
&=U^0(x,v)\;\delta_{p} (v-\frac{x}{\eps})
\end{aligned}
\end{equation}
with periodic boundary conditions on the torus in the $v$  variable.
Our objective is to analyze the $\eps \to 0$ limit of this problem and through this process
to calculate the effective equation satisfied by the weak limit of $u_{\eps}$. We note
that this is a hydrodynamic limit problem in the diffusive scaling for the kinetic
equation \eqref{maindiff}.

We prove.

\begin{theorem}
\label{enhdiffthm}
Under hypothesis \eqref{vecf} and \eqref{diffdata} we have the following 
asymptotic behavior
for $f_{\eps}$ as $\eps \to 0$:
\begin{align}
f_{\eps} (t,x,v) &\rightharpoonup u(t,x) \quad \text{weak-$\star$ in  
$L^\infty([0,\ \infty) ,\ L^2(\R^d, M_{p})$}
\\
f_{\eps} (t,x,v) &= u_{\eps}(t,x) + \eps g(t,x,v) + o(\eps) \, , \quad \text{in $\cD'$}
\end{align}
where $u$ of class \eqref{unifdiff} and  
$g \in L^{2} ( (0, \infty) \times \Omega ; H^{1} (\T^{d}))$ satisfy
respectively
\begin{align}
\del_{t} u -  \alpha \triangle_{x} u + 
\nabla_{x} \cdot \int_{\T^{d}} a(x,v) g (t,x,v) = 0
\label{basic}
\\
\alpha \triangle_{v} g  - \nabla_{v} \cdot a g  =
a \cdot \nabla_{x} u
\label{gsol}
\end{align}
The weak limit $u$ satisfies the effective diffusion equation
$$
\begin{aligned}
\del_{t} u &=  \alpha \sum_{i, j} \del_{x_{i}}  \left (
\Big  (\delta_{i j} +  \int_{\T^{d}}  \nabla_{v} \chi_{i} \cdot \nabla_{v} \chi_{j} dv  
\Big )
\del_{x_{j}} u \right )
\\
u (0,x) &= \int_{\T^{d}} U^{0} (x,v) \, dv
\end{aligned}
$$
where $\chi_{k}$, $k = 1, ... , d$, is the solution of the cell problem
\begin{equation}
\label{cellp}
\alpha \triangle_{v} \chi_{k} - \nabla_{v} \cdot a \chi_{k} = a \cdot e_{k} 
\end{equation}
\end{theorem}


\begin{proof} Let $f_{\eps}$  be defined as in \eqref{diffdsdef}. Then  $f_{\eps}$ satisfies 
the problem  \eqref{maindiff} and  $u_{\eps} = \int_{\T^{d}} f_{\eps}$.  
The proof is split in three steps:

\medskip
{\em Step 1~: Characterization of the weak limit.} 
 From \eqref{unifdiff} and Lemma \ref{lem1} we obtain uniform bounds
for $f_{\eps}$:
\begin{equation}
\label{fbounds}
\begin{aligned}
f_\ep &\in_{b} L^\infty([0,\ \infty) \,  ; \, L^2(\R^d, M_{p}) \, ,
\\
\big ( \nabla_{x} + \frac{1}{\eps} \nabla_{v}  \big ) f_{\eps}
&\in_{b} L^{2} ( (0, \infty) \times \R^{d} \, ; \ M_{p})
\end{aligned}
\end{equation}
Using (a slight variant of) Proposition \ref{prop1} we see that, along a subsequence if necessary,
$f_{\eps}$ satisfies
\begin{equation}
\label{wl}
\begin{aligned}
f_{\eps} &\rightharpoonup f \quad \text{weak-$\star$ in  
$L^\infty([0,\ \infty) ;\ L^2(\R^d , M_{p})$}
\\
f &\in L^\infty([0,\ \infty) \, ; \,  L^2(\R^d \times \T^{d})) \, ,
\end{aligned}
\end{equation}
that is
$$
\int \theta(t,x,v) df_{\eps}(t,x,v) \to \int \theta(t,x,v) f(t,x,v) dt dx dv
$$
for $\theta \in L^{1} ((0,\infty) ; L^{2} (\R^{d}, C_{p}))$.

Passing to the limit $\eps \to 0$ in \eqref{maindiff} and using \eqref{fbounds}
we see that $f$ satisfies
$$
\alpha \triangle_{v} f - \nabla_{v} \cdot a f = 0 \, , 
\quad \text{in $\cD'$},
$$
and for any test function $\theta$ 
$$
\begin{aligned}
<\nabla_{v} f , \theta > = \lim_{\eps \to 0} <\nabla_{v} f_{\eps}, \theta>
= 0
\, .
\end{aligned}
$$
Hence, $\nabla_{v}f = 0$ in $\cD'$ and 
\begin{equation}
\label{cl}
f(t,x,v) = \int_{\T^{d}} f dv = : u(t,x)
\end{equation}

\medskip
{\em Step 2~: Asymptotics of $f_{\eps}$.} Define next 
\begin{equation}
g_{\eps} (t,x,v) = \frac{1}{\eps} \big ( f_{\eps} (t,x,v)  - u_{\eps} (t,x) \big )
\, , \quad t \in \R_{+} \, ,  \, x \in \R^{d} \, , \, v \in \T^{d} \, ,
\end{equation}
where $u_{\eps} = \int_{\T^{d}} f_{\eps}$. We proceed along the lines
of Proposition \ref{prop3} replacing the bounds of that proposition by the bound 
\eqref{unifdiff} and accounting for the extra dependence in time. After minor
modifications in the proof we obtain for any $T>0$
\begin{equation}
\label{proptimeasym}
\begin{aligned}
g_{\eps} &\in_{b}  L^{2} \big ( (0,T) \, ; \, H^{-1} (\Omega , M_{p}) \big )
 \\
 g_{\eps} &\rightharpoonup g \quad 
 \text{ weak-$\star$ in $L^{2} \big ( (0,T) \, ; \, H^{-1} (\Omega , M_{p}) \big )$}
 \\
 g &\in L^{2} ( (0, \infty) \times \Omega \, ; \, H^{1} (\T^{d})) \, , 
 \quad \int_{\T^{d}} g = 0 \, .
\end{aligned}
\end{equation}
Accordingly, $f_{\eps}$ enjoys  the asymptotic expansion 
$$
f_{\eps} = u_{\eps} + \eps g + o(\eps) \quad \text{in $\cD'$} \, .
$$

On the other hand, on account of \eqref{enhdiff}, \eqref{maindiff} and \eqref{vecf},
it follows that  $u_{\eps}$ and $g_{\eps}$ satisfy
\begin{equation}
\begin{aligned}
\del_{t} u_{\eps} -  \alpha \triangle_{x} u_{\eps} + 
\nabla_{x} \cdot \int_{\T^{d}} a(x,v) \frac{f_{\eps } - u_{\eps} }{\eps} = 0
\end{aligned}
\end{equation}
and
\begin{equation}
\begin{aligned}
\alpha \triangle_{v} g_{\eps}  - \nabla_{v} \cdot a g_{\eps} =
\eps (\del_{t} f_{\eps} -  \alpha \triangle_{x} f_{\eps}) -
2 \alpha \nabla_{x} \cdot \nabla_{v} f_{\eps}
+ \nabla_{x} \cdot a f_{\eps}
\end{aligned}
\end{equation}
Using \eqref{wl}, \eqref{cl} and \eqref{proptimeasym}, we pass to the limit $\eps \to 0$
and deduce that $u$, $g$ satisfy \eqref{basic} and \eqref{gsol} respectively.

\medskip
{\em Step 3~: Characterization of the limit problem.} Due to its regularity the solution $g$ 
of \eqref{gsol} is unique and can be expressed in the form
$$
g = \nabla_{x} u(t,x) \cdot \chi (x,v)
$$
where $\chi = ( \chi_{1}, ... , \chi_{d})$ is the solution of the cell problem
\begin{equation}
\alpha \triangle_{v} \chi_{k} - \nabla_{v} \cdot a \chi_{k} = a \cdot e_{k} = a_{k}
\, .
\end{equation}
A direct computation shows that solutions of \eqref{cellp} satisfy the property
$$
\frac{1}{2} \int_{\T^{d}} (a_{i} \chi_{k} + a_{k} \chi_{i} ) dv =
- \alpha \int_{\T^{d}} \nabla_{v} \chi_{i} \cdot \nabla_{v} \chi_{k} dv
$$
and  \eqref{basic} may be written in the 
equivalent forms
$$
\begin{aligned}
\del_{t} u &= \sum_{i, j} \del_{x_{i}} \left (
\big  (\alpha \delta_{i j} - \frac{1}{2} \int_{\T^{d}} (a_{i} \chi_{j} + a_{j} \chi_{i})  \, dv   
\big )
\del_{x_{j}} u \right )
\\
&= \alpha \sum_{i, j} \del_{x_{i}} \left (
\Big  ( \delta_{i j} +  \int_{\T^{d}}  \nabla_{v} \chi_{i} \cdot \nabla_{v} \chi_{j} dv  
\Big ) \del_{x_{j}} u \right )
\end{aligned}
$$
The latter is a diffusion equation with positive definite diffusion matrix 
$$
\begin{aligned}
D_{i j} &=  \delta_{i j} +  \int_{\T^{d}}  \nabla_{v} \chi_{i} \cdot \nabla_{v} \chi_{j} dv
\\
\sum_{i j} D_{i j} \nu_{i} \nu_{j} &=
|\nu|^{2} + \int_{\T^{d}} | \nabla_{v} (\chi \cdot \nu) |^{2} dv \, ,  \quad 
\text{$\nu \in \R^{d}$},
\end{aligned}
$$
determined through the solution of \eqref{cellp}.
\end{proof}

%

\section{Appendix I}
\label{appendix}

We prove in the appendix a lemma that is used in the justification of multiscale
decompositions.

\begin{lemma}
\label{technical}
Let $\Omega$ be an open subset of $\R^{d}$, 
$\theta \in C_{c}(\Omega)$, $\varphi \in C(\T^{d})$, $\psi \in C(\T^{d})$,
and suppose that $\delta = \delta (\eps) \to 0$ as $\eps \to 0$. Then, as $\eps \to 0$,
\begin{align}
\int_{\Omega} \theta(x) \varphi(\frac{x}{\eps}) dx &\to
\int_{\T^{d}}\varphi(z)dz \; \int_{\Omega} \theta (x) dx
\label{techn1}
\\
\int_{\Omega} \theta(x) \varphi(\frac{x}{\eps}) \psi (\frac{x}{\eps \delta}) dx &\to
\int_{\T^{d}}\psi(y)dy \; \int_{\T^{d}}\varphi(z)dz \; \int_{\Omega} \theta (x) dx
\label{techn2}
\end{align}
\end{lemma}

\begin{proof}
Fix $\theta \in C_{c}(\Omega)$ and let
$\eps < \frac{1}{\sqrt{d}} \, \dist ( \supp \theta , \partial \Omega)$. 
We consider a cover of the support of the function $\theta$ 
by cubes $C_{k}$ centered at points $\chi_{k} \in \eps \Z^{d}$ of latteral size $\eps$.
The number of the cubes covering $\supp \theta$ satisfies $N \eps^{d} = O(|\supp \theta|)$,
and since  $\eps \sqrt{d} < \dist ( \supp \theta , \partial \Omega)$ the covering
can be arranged so that
$\supp \theta \subset \cup_{k=1}^{N} C_{k} \subset \Omega$ while 
$\frac{1}{\eps} \chi_{k} \in \Z^{d}$. We have
$$
\begin{aligned}
\int_{\Omega} \theta(x) \varphi(\frac{x}{\eps}) dx
& = \sum_{k=1}^{N} \int_{C_{k}} \theta(x) \varphi(\frac{x}{\eps}) dx
\\
&= \sum_{k=1}^{N} \eps^{d} \int_{\T^{d}} \theta (\chi_{k} + \eps z) 
\varphi ( \frac{1}{\eps} \chi_{k} + z) \, dz
\\
&= \sum_{k=1}^{N} \eps^{d} \int_{\T^{d}} \theta (\chi_{k} + \eps z) 
\varphi (  z) \, dz
\end{aligned}
$$
and thus
$$
\begin{aligned}
I &= \int_{\Omega} \theta(x) \varphi(\frac{x}{\eps}) dx -
\int_{\T^{d}}\varphi(y)dy \; \int_{\Omega} \theta (x) dx
\\
&= \sum_{k=1}^{N} \eps^{d} \int_{\T^{d}} \theta(\chi_{k} + \eps y) \varphi(y) dy
-
\int_{\T^{d}} \varphi(y) dy \sum_{k=1}^{N} \eps^{d} 
\int_{\T^{d}} \theta(\chi_{k} + \eps z) dz
\\
&= \int_{\T^{d}} \varphi(y) \sum_{k=1}^{N} \eps^{d} \int_{\T^{d}}
\left (  \theta(\chi_{k} + \eps y) - \theta(\chi_{k} + \eps z) \right ) dz \, dy
\end{aligned}
$$
Using the uniform continuity of $\theta$ and that $N = O(\frac{1}{\eps^{d}})$ we 
deduce $I \to 0$ as $\eps \to 0$ and \eqref{techn1}.

Next observe that
$$
\begin{aligned}
&\int_{\Omega} \theta(x) \varphi(\frac{x}{\eps}) \psi (\frac{x}{\eps \delta}) dx -
\int_{\T^{d}}\psi(y)dy \; \int_{\T^{d}}\varphi(z)dz \; \int_{\Omega} \theta (x) dx
\\
&\quad =
\int_{\Omega} \theta(x) \varphi(\frac{x}{\eps}) \psi (\frac{x}{\eps \delta}) dx
- \int_{\T^{d}} \psi(y) dy \; \int_{\Omega} \theta(x) \varphi(\frac{x}{\eps}) dx
\\
&\qquad + \int_{\T^{d}} \psi(y)dy 
\left(  
\int_{\Omega} \theta(x) \varphi(\frac{x}{\eps}) dx - \int_{\T^{d}}\varphi(z) dz
\; \int_{\Omega} \theta(x)dx
\right)
\\
= I_{1} + I_{2}
\end{aligned}
$$
and that $I_{2} \to 0$ as $\eps \to 0$.

Consider now a covering of $\supp \theta$ by cubes $\bar C_{k}$ centered at points
$\chi_{k} \in \eps \delta \Z^{d}$ of latteral size $\eps \delta$. 
As in the preceding argument we can arrange the cubes so that
$\supp \theta \subset \cup_{k=1}^{\bar N} \bar C_{k} \subset \Omega$ and their number 
$\bar N$ satisfies 
 $\bar N (\eps \delta)^{d} = O(|\supp \theta|)$. We have
$$
\begin{aligned}
\int_{\Omega} \theta(x) \varphi(\frac{x}{\eps}) \psi (\frac{x}{\eps \delta}) dx
&=
\sum_{k=1}^{{\bar N}} \int_{\bar C_{k}} 
\theta(x) \varphi(\frac{x}{\eps}) \psi (\frac{x}{\eps \delta}) dx
\\
&=
\sum_{k=1}^{{\bar N}} (\eps \delta)^{d} \int_{\T^{d}} 
\theta(\bar \chi_{k} + \eps \delta z) 
\varphi(\frac{1}{\eps} \bar \chi_{k} + \delta z) \psi (z) dz
\end{aligned}
$$
and 
$$
\begin{aligned}
I_{1} &=
\sum_{k=1}^{{\bar N}} \int_{\bar C_{k}} 
\theta(x) \varphi(\frac{x}{\eps}) \psi (\frac{x}{\eps \delta}) dx
- \int_{\T^{d}} \psi(z) dz \; \sum_{k=1}^{\bar N}\int_{\bar C_{k}} 
\theta (x) \varphi(\frac{x}{\eps}) dx
\\
&=
\int_{\T^{d}} \psi(z) \sum_{k=1}^{\bar N} (\eps \delta)^{d} \int_{\T^{d}}
\Big ( \theta (\bar \chi_{k} +\eps \delta z) 
\varphi( \frac{1}{\eps}\bar \chi_{k} + \delta z)
\\
&\qquad \qquad \qquad \qquad  \qquad \qquad -
\theta (\bar \chi_{k} +\eps \delta y) 
\varphi( \frac{1}{\eps}\bar \chi_{k} + \delta y)
\Big ) \, dy \, dz
\end{aligned}
$$
Again, since $\bar N = O(\frac{1}{(\eps \delta)^{d}})$ and $\lim_{\eps\to 0}\delta = 0$, 
we deduce $I_{2} \to 0$ as $\eps \to 0$
and \eqref{techn2}.
\end{proof}

\begin{remark} Relation \eqref{techn1} is classical (see \cite{BLL}) and is only proved here
as a precursor to the proof of \eqref{techn2}. Equation \eqref{techn2} indicates that oscillations
of different scales do not correlate and suggests that the definition \eqref{defn3sc} 
is meaningful.

Both equations can be extended for test functions $\theta \in C(\bar \Omega)$ provided
that $\Omega \subset \R^{d}$ is a bounded open set and its boundary $\partial \Omega$  
is of finite $d-1$  Hausdorff dimension.
\end{remark}

As an application we prove \eqref{ex1}.

\begin{lemma}
\label{lemapen}
We have
$$
\delta_p \big ( v - \frac{x}{\eps} \big ) \rightharpoonup 1
\qquad \text{weak-$\star$ in $L^\infty(\Omega, M_p)$ }
$$
\end{lemma}

\begin{proof}
We need to show that for $\theta \in L^1(\Omega, C_p)$ we have
$$
\int_{\Omega} \theta (x, \frac{x}{\eps}) dx \to \int_\Omega \int_{\T^d}
\theta(x,v) dx dv
$$
Equation \eqref{techn1} justifies that for $\theta = \chi(x) \otimes \varphi(v)$
a tensor product with $\chi \in C_c(\Omega)$ and $\varphi \in C_p$ and
by a density argument also for $\chi \in L^1 (\Omega)$, $\varphi \in C_p$.

To complete the proof we need to show that finite sums of tensor products 
$\sum_j \chi_j \otimes \varphi_j$ are dense in $L^1 (\Omega, C_p)$. Fix 
$\theta \in L^1 (\Omega, C_p)$ and consider a decomposition of the torus
 $\T^d$ into squares of size $1/n$. Take a partition of unity 
$\varphi_i \in C_p$, $i = 1, ..., n^d$, with each $\varphi_i$ supported
in a square of size $2/n$ and $\sum_i \varphi_i =1$. Let $v_i$ be the
center of each square and define
$$
\theta_n (x,v) = \sum_i \theta (x, v_i) \varphi_i (v)
$$
Clearly, $\theta_n$ is a sum of tensor products. Now define
$$
\sup_{v \in \T^d} |\theta(x,v) - \theta_n (x,v) | 
\le \sup_{v \in \T^d , |h| < \frac{2}{n}} |\theta(x,v) - \theta (x,v+h) |
=: g_n(x)
$$
and thus
$$
\| \theta - \theta_n \|_{L^1(C_p)} \le \int_\Omega |g_n(x)| dx
$$
Note that $g_n(x) \to 0$ for a.e. $x\in \Omega$ and that 
$|g_n(x)| \le 2 \sup_{v\in \T^d} |\theta(x,v)|$. The latter
is an $L^1$ function by the very definition of $\theta$, and the
dominated convergence theorem implies $\int_\Omega |g_n|dx \to 0$. 
\end{proof}

%

\section{Appendix II: Some basic results of ergodic theory}
The purpose of this appendix is to recall some well known properties
of the classical ergodic theory
for the projection on the kernel
\[
K=\{f\in L^2(\T^d)\,|\;a(v)\cdot\nabla_v f(v)=0\},
\]
where the last equality is of course in the sense of distributions.

Let us define the characteristics associated with $a$ which are the
solutions on $\T^d$ of the following differential equation
\[
\partial_t T(t,v)=a(T(t,v)),\quad T(0,v)=v.
\]
Then assuming that
\begin{equation}
a\in W^{1,\infty}(\T^d),\quad \nabla_v\cdot a=0,\label{hypa}
\end{equation}
the characteristics are well defined and for a fixed $t$, the
transform $v\rightarrow T(t,v)$ is a mesure preserving homeomorphism
of $\T^d$. We then have the well-known theorem (see Sinai \cite{Si}
for more details)
\begin{theorem}
For every $f\in L^p(\T^d)$ with $1\leq p<\infty$, there exists a
unique function in $L^p(\T^d)$, denoted by $Pf$, such that 
\[
\frac{1}{t}\int_0^t f(T(s,v))\,ds\longrightarrow Pf(v),\ \hbox{as}\
t\rightarrow \infty,\ 
\hbox{strongly\ in}\ L^p(\T^d).
\]
Moreover $Pf$ satisfies in the sense of distribution
\[
a(v)\cdot\nabla_v Pf(v)=0,
\]
and if $f\in L^2(\T^d)$, then $Pf$ is exactly the orthogonal
projection of $f$ on $K$.
\label{ergodic}
\end{theorem}

This immediately implies the

\begin{corollary}
The orthogonal projection on $K$ may be extended as an operator on
$L^p(\T^d)$ for every $1\leq p<\infty$. In addition if $f\in L^2\cap
L^p(\T^d)$ with $1\leq p\leq\infty$ ($p=\infty$ allowed), then $P_K
f$ also belongs to $L^2\cap
L^p(\T^d)$.
\end{corollary}
{\bf Proof of Theorem \ref{ergodic}.} This proof exactly corresponds
to the one in \cite{Si} in the particular case which we consider.

Notice that if, for $f\in L^p(\T^d)$, $\frac{1}{t}\int_0^t
f(T(s,v))\,ds$ converges to $Pf$ then trivially
\[
Pf(T(t,v))=Pf(v)\quad \forall\ t.
\]
Therefore we automatically have in the sense of distribution that
\[
a(v)\cdot\nabla_v Pf=0.
\]
Take now $f$ in $L^p$ and assume
first that there exists $g\in L^p$ with $a\cdot\nabla_v g=0$ and 
$h\in W^{1,p}(\T^d)$ such that
\[
f=g+a\cdot\nabla_v h.
\]
Then notice that in the sense of distribution
\[
\partial_t\left(g(T(t,v))\right)=a(T(t,v))\cdot\nabla_v g(T(t,v))=0,
\]
and so
\[
g(T(t,v))=g(T(0,v))=g(v).
\]
On the other hand
\[
a(T(t,v))\cdot\nabla_v h(T(t,v))=\partial_t\left(h(T(t,v))\right),
\]
and therefore
\[
\frac{1}{t}\int_0^t f(T(s,v))\,ds=g(v)+\frac{h(T(t,v))-h(v)}{t}.
\]
Consequently in this case $\frac{1}{t}\int_0^t f(T(s,v))\,ds$
converges to $g$ which is unique as a consequence. This proves the
theorem on the set
\[
L_p=\{g+a(v)\cdot\nabla_v h(v)\,|\;h\in W^{1,\infty}(\T^d),\ g\in
L^p(\T^d)\ \hbox{with}\ a\cdot\nabla_v g=0\}.
\]
Let us first prove that $L_2$ is dense in $L^2(\T^d)$.
If $L_2$ is not dense, then there exists $f\in
L^2\setminus \{0\}$ orthogonal to $L_2$. This implies that for all
$h\in W^{1,\infty}(\T^d)$
\[
\int_{\T^d} f(v)\,a(v)\cdot\nabla_v h(v)\,dv=0, 
\]
or in other words $f\in K$. But $K\subset L_2$ and $f$ should
consequently be
orthogonal to $K$, which is impossible. Notice that $Pf$ necessarily is the orthogonal
projection on $K$ as $a\cdot\nabla h$ belongs to $K^\perp$.

Now for any $f\in L^p$. If $1\leq p<2$, take $g_n+a\cdot\nabla_v h_n=f_n\in L_2$ converging
toward $f$ in $L^p$ (first take ${\hat f}_n \in L^2$ 
and then select $f_n$ by diagonal extraction). We have that
\[
\begin{aligned}
\|\frac{1}{t}\int_0^t f(T(s,v))\,ds &- \frac{1}{t'}\int_0^{t'
}f(T(s,v))\,ds\|_{L^p}
\\
&\leq
2\|f-f_n\|_{L^p}
+\left(\frac{1}{t}+\frac{1}{t'}\right)\;\|h_n\|_{L^p}.
\end{aligned}
\]
So the sequence $\frac{1}{t}\int_0^t f(T(s,v))\,ds$ is of Cauchy in
$L^p$ and hence converges to a unique limit $Pf$. 

Finally if $f\in L^\infty(\T^d)$, then $f\in L^2(\T^d)$ and 
$\frac{1}{t}\int_0^t f(T(s,v))\,ds$
converges to $Pf$ in $L^2$. As the first quantity is uniformly bounded
in $L^\infty$, $Pf\in L^\infty$ and the convergence holds in every
$L^p$, $p<\infty$. By interpolation, one eventually obtains the
desired result for $f\in L^p(\T^d)$.\cqfd

\section*{Acknowledgements}
PEJ was partially supported by the HYKE european network. 
AET is partially supported by the National Science Foundation.



\begin{thebibliography}{10}

\bibitem{Allaire} Allaire, G., Homogenization and two-scale convergence, 
{\em SIAM J. Math. Anal.} {\bf 23} (1992), 1482-1518.

\bibitem{AlVa} Allaire, G., Vanninathan, M., 
Homogenization of the Schr\"odinger equation with a time oscillating
potential.  
{\em Discrete Contin. Dyn. Syst. Ser. B}  {\bf 6}  (2006),  no. 1, 1--16 

\bibitem{AHZ} Amirat, Y., Hamdache, K., Ziani, A., 
Homog\'en\'eisation d'\'equations hyperboliques du premier ordre et 
application aux \'ecoulements miscibles en milieu poreux. 
(French) [Homogenization of a system of first-order hyperbolic
  equations 
and application to miscible flows in a porous medium]  
{\em Ann. Inst. H. Poincar\'e Anal. Non Linéaire}  {\bf 6}  
(1989),  no. 5, 397--417.

\bibitem{AM} Avellaneda, M. and A.J. Majda,
An integral representation and bounds on the effective diffusivity in passive advection
by laminar and turbulent flows.
{\em Comm. Math. Physics} {\bf 138} (1991), 339-391.

\bibitem{BLL} Bensoussan, A., J.L. Lions and G. Papanicolaou, {\em Asymptotic Analysis
for Periodic Structures}, North Holland, Amsterdam, 1978.

\bibitem{Brenier} Brenier,  Y., Remarks on some linear hyperbolic
  equations with oscillatory coefficients.   {\em Third International
  Conference on Hyperbolic Problems, Vol. I, II (Uppsala, 1990)},
  119--130, {\em Studentlitteratur, Lund}, 1991.

\bibitem{Capdeboscq1} Capdeboscq, Y., Homogenization of a spectral
  problem with drift, {\em Proc. Roy. Soc. Edinburgh Sect. A}  {\bf 132}
  (2002),  no. 3, 567--594.  

\bibitem{Capdeboscq2} Capdeboscq, Y., {\it Homog\'en\'eisation des
  mod\`eles de diffusion en neutronique}, Th\`ese Universit\'e Paris
  6, 1999.

\bibitem{E} E, W., 
Homogenization of linear and nonlinear transport equations, 
{\em Comm. Pure Appl. Math.} {\bf 45} (1992), 301-326.

\bibitem{FP} Fannjiang, A. and G. Papanicolaou,
Convection enhanced diffusion for periodic flows,
{\em SIAM J. Appl. Math.} {\bf 54} (1994), 333-408.

\bibitem{GMMP} G\'erard, P., Markowich P. A., Mauser, N. J., 
Poupaud, F., Homogenization limits and Wigner transforms.  
{\em Comm. Pure Appl. Math.}  {\bf 50}  (1997),  no. 4, 323--379. 

\bibitem{GouPou} Goudon, T. and Poupaud, F., 
Homogenization of transport equations: weak mean field approximation.  
{\em SIAM J. Math. Anal.}  {\bf 36}  (2004/05),  no. 3, 856--881

\bibitem{Hamdache} Hamdache, K., { Homog\'en\'eisation non locale
  d'\'equations hyperboliques}, in {\it Non linear pde's and their
  applications, Coll\`ege de France seminar}, vol. XII,Pitman
  Res. Notes in Math., {\bf 302}, (Longman Sci. Tech., Harlow, 1994)
  pp 97--112.
 
\bibitem{Hormander} Hormander, L., {\em The Analysis of Linear Partial Differential Operators}.
Springer, New York, 1990.

\bibitem{HX} Hou, T.Y. and X. Xin,
Homogenization of linear transport equations with oscillatory vector fields,
{\em SIAM J. Appl. Math.} {\bf 52} (1992), 34-45.

\bibitem{MPP} McLaughlin D.W., G.C. Papanicolaou and O.R. Pironneau,
Convection of microstructures and related problems,
{\em SIAM J. Appl. Math.} {\bf 45} (1985), 780-797.

\bibitem{Nguetseng} Nguetseng G., A general convergence result for a 
functional related to the theory of homogenization. {\em SIAM J. Math. Anal.}
{\bf 3}, 1989, 608--623.

\bibitem{Si} Sinai, Ya. G., Dynamical Systems III, {\em Springer
  Verlag}, New-York Heidelberg Berlin, 1989.

\bibitem{Schwartz} Schwartz, L., {\em Th\'eorie des Distributions}.
Hermann, Paris, 1966.

\bibitem{Tartar1} Tartar, L., {Remarks on homogenization} in {\it
  Homogenization and effective moduli of material and media}, IMA
  Vol. in Math. and Appl. (Springer, 1986) pp. 228--246.

\bibitem{Tartar2} Tartar L., {Nonlocal effects induced by
  homogenization}, in {\it Partial differential equations and the
  calculus of variations, Essays of Mathematical analysis in honor of
  E. De Giorgi}, Vol. II, Progr. Nonlinear Differential Equations
  Appl.,{\bf 2} (Birkhauser, 1989) pp. 925--938.

\end{thebibliography}
\end{document}